\documentclass[reqno,10pt]{amsart}
\usepackage{CJK}
\usepackage{hyperref}
\usepackage{cite}
\usepackage{amsmath}
\usepackage{mathrsfs}
\usepackage{extarrows}
\usepackage{indentfirst}
\usepackage{color}

\makeatletter
\@namedef{subjclassname@2020}{%
	\textup{2020} Mathematics Subject Classification}
\makeatother

\numberwithin{equation}{section}

\newtheorem{theorem}{Theorem}[section]
\newtheorem{lemma}[theorem]{Lemma}
\newtheorem{definition}[theorem]{Definition}

\newtheorem{remark}[theorem]{Remark}
\newtheorem{corollary}[theorem]{Corollary}

\allowdisplaybreaks

\newcommand{ \mint }{ {\int\hspace{-0.38cm}-}}

\begin{document}
	
	\title[\hfil Higher Sobolev regularity for mixed problems] {Higher Sobolev regularity for mixed local and nonlocal equations with nonstandard growth}
	
		\author[Y. Fang, D. Li and C. Zhang  \hfil \hfilneg]
	{Yuzhou Fang, Dingding Li  and Chao Zhang$^*$}
	
	\thanks{$^*$Corresponding author.}
	
	\address{Yuzhou Fang  \hfill\break School of Mathematics, Harbin Institute of Technology, Harbin 150001, China}
	\email{18b912036@hit.edu.cn}
	
	\address{Dingding Li \hfill\break School of Mathematics, Harbin Institute of Technology, Harbin 150001, China}
	\email{a87076322@163.com}

	\address{Chao Zhang  \hfill\break School of Mathematics and Institute for Advanced Study in Mathematics, Harbin Institute of Technology, Harbin 150001, China}
	\email{czhangmath@hit.edu.cn}

	\subjclass[2020]{35B65, 35D30, 35J60, 35R11}
	\keywords{Sobolev regularity; mixed local and nonlocal $p$-Laplacian; finite difference quotients}
	
\begin{abstract}
We develop a systematic study of the interior Sobolev regularity of weak solutions to the mixed local and nonlocal $(p,q)$-Laplace equations, which model the superposition of two stochastic processes such as Brownian motion
and L\'{e}vy flight. To be precise, we show that the weak solution $u$ belongs to $W^{2, 2}_{\rm loc}$ Sobolev spaces in the subquadratic case, while $|\nabla u|^{\frac{p-2}{2}}\nabla u$ is of the class $W^{1, 2}_\mathrm{loc}$ in the superquadratic scenario, both of which coincide with that of the classical $p$-Laplace equations. Moreover, an improved higher fractional differentiability and integrability result $u\in W^{1+\beta, \gamma}_\mathrm{loc}$ is proved in the full range $p\in (1, \infty)$ for any $\gamma\in [\max\{p, 2\}, \infty)$ and $\beta\in(0, \frac 2\gamma)$. The main analytical tools are the finite-diﬀerence quotients, carefully tailored energy estimates, and
nonlocal tail estimates. As far as we know, our results are new within the context of such mixed problems, even in the $p$-growth setting.
\end{abstract}

	\maketitle

\section{Introduction}
\thispagestyle{empty}
\label{sec1}

In this paper, we study the higher Sobolev regularity, in both differentiability and integrability, of weak solutions to the mixed local and nonlocal equation
\begin{align}
	\label{main}
	-\Delta_p u+(-\Delta_q)^s u=0\quad\text{in }\Omega,
\end{align}
where
$$
s\in (0,1), \quad sq<p\le q<\infty
$$
and $\Omega\subset\mathbb{R}^n$, $n\ge2$, is a bounded domain. The local $p$-Laplace operator $\Delta_p$ and the fractional $q$-Laplace operator $(-\Delta_q)^s$ are defined, respectively, by
\begin{align*}
	\Delta_pu=\mathrm{div}\left( |\nabla u|^{p-2}\nabla u\right)
\end{align*}
and
\begin{align*}
	(-\Delta_q)^su(x)=2 \mathrm{P.V.}\int_{\mathbb{R}^n}\frac{|u(x)-u(y)|^{q-2}(u(x)-u(y))}{|x-y|^{n+sq}}\,dy.
\end{align*}
Here $\nabla u$ denotes the gradient of $u$ and $\mathrm{P.V.}$ denotes the Cauchy principal value. Such equations, which simultaneously exhibit the classical and fractional features, are not scale-invariant and therefore cannot be regarded as a simple superposition.

In recent years, mixed local and nonlocal problems have received increasing attention owing to their far-ranging applications in the real
world. On one hand, the analysis of mixed operators is motivated by the superposition of two stochastic processes such as Brownian motion
and L\'{e}vy flight. On the other hand, such problems arise in optimal search strategies, animal foraging, and biomathematics, where they describe diffusion patterns across multiple scales; see \cite{DV21,DLV24} and the references therein.

\subsection{Overview of related results}

The prototype $p$-Laplace equation $-\Delta_pu=0$ has been studied extensively, and the regularity of the gradient of its weak solutions is now well understood. We recall several representative local higher-regularity results for a $p$-harmonic function $u\in W^{1,p}$:
\begin{itemize}
	\item $\nabla u\in C^{0, \alpha}_{\rm loc}$ for $1<p<\infty$ (see e.g. \cite{LJ83, T84, U77, U68});
	
	\smallskip
	
	\item $\nabla u\in W^{1, p}_{\rm loc}$ for $1<p\le2$ in \cite[Chapter 4]{L19} and $|\nabla u|^{\frac{p-2}{2}}\nabla u\in W^{1,2}_{\rm loc}$ for $2\le
p<\infty$ in \cite{BI87};
	
	\smallskip
	
	\item $\nabla u\in W^{1, 2}_{\rm loc}$ under the Cordes condition $1<p<3+\frac{2}{n-2}$ by Manfredi and Weitsman \cite{MW88} (see \cite{HS22} for a global version of this);
	
	\smallskip
	
	\item $\nabla u\in \mathcal{N}^{\frac{2}{p},p}$ (a Nikol'skii space) for $p\ge2$ by Mingione \cite{M03}; see also \cite{PAMS18} for the improved result $\nabla u\in \mathcal{N}^{\theta,\frac{2}{\theta}}$ with any $\theta\in \left[ \frac{2}{p},\frac{2}{p-1}\right) $ in the case $p\ge3$.
\end{itemize}
Several of these results extend to more general second-order operators, systems, or equations with nontrivial source terms. For the nonhomogeneous equation $-\Delta_pu=f$, Cianchi and Maz'ya \cite{CM18} proved, under minimal regularity assumptions on the domain, that $|\nabla u|^{p-2}\nabla u\in W^{1,2}$ if and only if $f\in L^2$; see \cite{CM19} for the corresponding result for the $p$-Laplace system. Dong et al. \cite{DPZZ20} subsequently refined the nonlinear gradient expressions through fundamental inequalities. They proved that $|\nabla u|^{\frac{p-2+s}{2}}\nabla u\in W^{1,2}_\mathrm{loc}$ for $p$-harmonic functions whenever $s>\max\{-1-\frac{p-1}{n-1},2-p-\frac{n}{n-1}\}$; see also \cite{S22} for the range $s>-1-\frac{p-1}{n-1}$.

For the nonlocal analogue of the $p$-Laplace equation, namely $(-\Delta_p)^su=0$, or, more generally,
\begin{align}
	\label{fp}
	\mathrm{P.V.}\int_{\mathbb{R}^n}|u(x)-u(y)|^{p-2}(u(x)-u(y))\mathcal{K}(x,y)dy=:\mathcal{L}u(x)=0
\end{align}
where $\mathcal{K}(x,y)\approx|x-y|^{-n-sp}$, a broad qualitative and quantitative theory has developed over the past two decades, beginning with the pioneering linear results in \cite{S06,CS07}. For nonlinear operators with $p$-growth, Di Castro, Kuusi, and Palatucci \cite{DKP14,DKP16} extended the De Giorgi--Nash--Moser theory to \eqref{fp}, proving local boundedness, H\"older continuity, and a Harnack inequality. Cozzi \cite{CM17} obtained related results using fractional De Giorgi classes. These methods have since been applied to obstacle problems \cite{KKP16}, fractional superharmonic functions and semicontinuity representations \cite{KKP17}, equations with measure data \cite{KMS15}, and the Wiener criterion for nonlocal Dirichlet problems \cite{KLL23}.

Turning to higher regularity for fractional $p$-Laplace equations, Brasco and Lindgren \cite{BL17} used finite-difference quotients to establish higher differentiability of weak solutions in the superquadratic case $p\ge2$, under an additional differentiability assumption at infinity. In particular, for sufficiently large $s$, they proved that $u$ is differentiable and $\nabla u\in W^{\tau,p}_\mathrm{loc}$. Still for $p\ge2$, Brasco, Lindgren and Schikorra \cite{BLS18} combined difference quotients with a Moser-type iteration to obtain higher integrability and, through Morrey-type embeddings, improved H\"older regularity; see \cite{GL24} for the subquadratic case $1<p\le2$. Diening et al. \cite{DKLN23} established higher differentiability for inhomogeneous fractional $p$-Laplace equations under various assumptions on the source term. B\"{o}gelein et al. \cite{BDL2402} enlarged, for $p\ge2$, the range of $s$ for which weak derivatives of $(s,p)$-harmonic functions exist, thereby extending \cite{BL17}; the singular case is treated in \cite{BDL2401}. For $p=2$, Diening et al. \cite{DKLN24} proved gradient H\"older continuity together with gradient potential estimates. Further results on nonlocal $p$-Laplacian problems may be found in \cite{BK23,DN23,KNS22,N23,KMS15b} and the references therein.

For the linear combination of local and nonlocal Laplace operators, namely,
\begin{align}
	\label{mix}
	-\Delta u+(-\Delta)^su=0,
\end{align}
probabilistic and analytic methods have yielded Harnack inequalities for globally nonnegative solutions and H\"older continuity for \eqref{mix} and its parabolic counterpart; see \cite{BBCK09,CKSV12,F09}. More recently, purely analytic approaches have addressed existence and nonexistence, boundedness, $C^{1,\alpha}$ and Sobolev regularity, strong maximum principles, and symmetry properties; see, for instance, \cite{AC21,BVDV21,BDVV21,DLV24}.

In the nonlinear setting $p\neq2$, Garain and Kinnunen \cite{GK24} derived local regularity results for
\begin{align}
	\label{mp}
	-\Delta_p u+\mathcal{L}u=0,
\end{align}
including boundedness, Harnack estimates, H\"older continuity, and lower semicontinuity of weak supersolutions, by means of De Giorgi--Nash--Moser iteration. For mixed local and nonlocal problems with measure data, Byun and Song \cite{BS23} obtained Wolff-potential estimates for solutions.

De Filippis and Mingione \cite{DM24} initiated the study of interior gradient H\"older regularity for mixed local and nonlocal $p$-Laplacian problems. In fact, they considered the more general functional with $(p,q)$-growth
\begin{align}
	\label{pq}
	u\longmapsto \int_\Omega\left( |\nabla
u|^p-fu\right)\,dx+\int_{\mathbb{R}^n}\int_{\mathbb{R}^n}\frac{|u(x)-u(y)|^q}{|x-y|^{n+sq}}\,dxdy
\end{align}
where $p,q>1$, $s\in(0,1)$, and $p>sq$, and used a perturbative argument to prove that its minimizers are locally $C^{1,\alpha}$ and globally almost Lipschitz continuous. Inspired by the difference-quotient approach of \cite{BLS18}, Garain and Lindgren \cite{GL2401} established improved H\"older regularity in the superquadratic case $p=q\ge2$. Related regularity results for $sq\ge p$ appear in \cite{DFZ24}. In addition, Byun, Kumar and Lee \cite{BKL24} obtained an optimal Calder\'{o}n--Zygmund estimate for mixed local and nonlocal equations with nonstandard growth, represented by the prototype
\begin{align*}
	-\mathrm{div}\left( |\nabla u|^{p-2}\nabla u\right)+(-\Delta_q)^su=-\mathrm{div}\left( |F|^{p-2}F\right).
\end{align*}
For further results on mixed problems, we refer to \cite{AC21,CKS11,FSZ22,GK24,SZ22,CSYZ24} and the references therein.

\subsection{Main results}

For \eqref{main}, one could expect that the regularity of weak solutions to \eqref{main} goes beyond that allowed by purely nonlocal operators, since the local operator $\Delta_p$ plays the leading role in view of $sq<p$.
To this end, we implement a difference-quotient method to obtain higher Sobolev estimates in both differentiability and integrability. More precisely, we differentiate the equation at the discrete level using carefully chosen test functions. The main difficulty is to control both the local integral and the nonlocal tail generated by $(-\Delta_q)^s$ in a manner compatible with the estimates furnished by the $p$-Laplacian. This yields matching powers of the increment $\delta_hu$ and the step size $|h|$, allowing us to pass rigorously from finite differences to weak derivatives. Consequently, we obtain stronger Sobolev regularity than is generally available for the purely fractional equation $(-\Delta_q)^su=0$.

The main contributions of this paper are threefold. First, we derive estimates for the first-order difference quotients of $|\nabla u|^{\frac{p-2}{2}}\nabla u$ for $p\ge2$ and of $\nabla u$ for $1<p\le2$, respectively. These estimates further yield their $W^{1,2}$-regularity. However, this cannot be achieved for the purely fractional $q$-Laplace equation, which does display the regularizing impact of $\Delta_p$. Second, the relation between difference quotients and weak derivatives in $W^{1,p}$ yields an increment of order $|h|^p$, which is sharper than the order $|h|^{sq}$ supplied by the $W^{s,q}$-energy. As a consequence, the order of differentiability of solutions in Theorem \ref{th2} below is higher than that of purely fractional problems. On the other hand, compared to the Sobolev regularity of $\nabla u$ for purely (local) $p$-Laplacian mentioned above due to Mingione, Mi\'skiewicz in \cite{M03, PAMS18}, where $\nabla u$ is shown to be almost $W^{\frac{2}{p}, p}$-regular in \cite{M03} for $p\ge2$ and almost $W^{\theta, \frac{2}{\theta}}$-regular with any $\theta\in \left[ \frac{2}{p}, \frac{2}{p-1}\right) $ in \cite{PAMS18} for $p\ge3$, our results improve the gradient integrability and provide an analogous theory in the singular range. Finally, it is worth pointing out that even though the classical operator $\Delta_p$ plays the dominant role, we cannot directly apply the methods used for the $p$-Laplace equation, owing to the nonlocal effects and a lack of scale-invariance. Overcoming these obstacles is essential to establish the desired Sobolev regularity.

Throughout this paper, we denote
$$
\mathfrak{M}_{r}=\|u\|_{L^\infty(B_{r}(x_0))}+\mathrm{Tail}(u;x_0,r)
$$
and
\begin{align*}
\mathrm{Tail}(u;x_0,r)=\left( r^{sq}\int_{\mathbb{R}^n\backslash B_r(x_0)}\frac{|u|^{q-1}}{|x-x_0|^{n+sq}}\,dx\right)^\frac{1}{q-1}
\end{align*}
for any $r>0$.

We now state our main results. The first theorem gives $u\in W^{2,2}_{\mathrm{loc}}$ for $1<p\le2$ and $|\nabla u|^{\frac{p-2}{2}}\nabla u\in W^{1,2}_{\mathrm{loc}}$ for $p\ge2$, in agreement with the corresponding regularity for $-\Delta_pu=0$ established in \cite{MW88,BI87}.

\begin{theorem}
	\label{th1}
	Let $s\in (0,1)$ and $u\in W^{1,p}_\mathrm{loc}(\Omega)\cap L^{q-1}_{sq}(\mathbb{R}^n)$ be a weak solution of \eqref{main}.
	\begin{itemize}
		
		\item[(i)] If $p\ge2$ and $u$ is locally bounded, then $|\nabla u|^{\frac{p-2}{2}}\nabla u\in
W^{1,2}_{\mathrm{loc}}(\Omega,\mathbb{R}^n)$. Moreover, for any $B_{2R}(x_0)\subset\subset\Omega$, we have
		\begin{align*}
			\int_{B_{R/2}(x_0)}\left| \nabla\left( |\nabla u|^{\frac{p-2}{2}}\nabla u\right) \right|^2\,dx\le C\left(1+ \|\nabla u\|_{L^p(B_{2R}(x_0))}+\mathfrak{M}_{2R}\right)^q.
		\end{align*}

		\item[(ii)] If $1<p\le2$ and $u$ is locally Lipschitz continuous, then $u\in W^{2,2}_\mathrm{loc}(\Omega)$.
		Moreover, for any $B_{2R}(x_0)\subset\subset\Omega$, we have
		\begin{align*}
			[\nabla u]^2_{W^{1,2}\left(B_{R/8}(x_0)\right)}\le C\left( 1+\mathfrak{M}_{2R}+\|\nabla u\|_{L^\infty(B_{2R}(x_0))}\right)^{2-2p+2q}.
		\end{align*}
	\end{itemize}
	Here the positive constant $C$ depends only on $n, s, p, q$ and $R$.
\end{theorem}

The local Lipschitz assumption in the range $1<p\le2$ is not restrictive: this regularity was established by De Filippis and Mingione \cite{DM24}. The Ishii--Lions method also yields a bound for $\nabla u$ in terms of the local supremum of $u$ and its tail. In our argument, Lipschitz continuity gives the pointwise estimate $|u(x+h)-u(x)|\le C|h|$ and allows us to extract an additional power of $|h|$ from the fractional term. The resulting estimate for the first-order difference quotient of $\nabla u$ is therefore uniform in $h$.

For $1<p\le2$, a Moser-type iteration improves the natural regularity $u\in W^{1,p}$ to $u\in W^{1,\gamma}$ for every $\gamma\ge p$, as stated next. In the superquadratic case, the corresponding higher-integrability result follows from the argument of Garain and Lindgren \cite{GL2401}.

\begin{theorem}
	\label{th3.13}
	Let $p\in(1,2]$ and $s\in (0,1)$. Assume that $u\in W^{1,p}_\mathrm{loc}(\Omega)\cap L^{q-1}_{sq}(\mathbb{R}^n)$ is a locally bounded weak
solution of \eqref{main}. Then we have
	\begin{align*}
		u\in W^{1,\gamma}_\mathrm{loc}(\Omega)\quad\text{for any }\gamma\in [p,\infty).
	\end{align*}
\end{theorem}

Our third result establishes fractional differentiability of $\nabla u$: for $p\in(1,\infty)$, the gradient is locally almost $W^{2/\gamma,\gamma}$ for every admissible $\gamma$.

\begin{theorem}
	\label{th2}
	Let $s\in (0,1)$ and let $u\in W^{1,p}_\mathrm{loc}(\Omega)\cap L^{q-1}_{sq}(\mathbb{R}^n)$ be a weak solution of \eqref{main}.
	\begin{itemize}
		
		\item[(i)] Let $p\ge2$ and suppose $u$ is locally bounded and belongs to $W^{1,\gamma}_{\mathrm{loc}}(\Omega)$. Then
		\begin{align*}
			u\in W^{1+\alpha,\gamma}_{\rm loc}(\Omega)
		\end{align*}
		for any $\alpha\in \left( 0,\frac{2}{\gamma}\right)$. Moreover, for $B_R(x_0)\subset\subset\Omega$,
		\begin{align*}
			[ \nabla u]^\gamma_{W^{\alpha,\gamma}(B_{R/4}(x_0))}\le C\left(1+\mathfrak{M}_R+\|\nabla u\|_{L^\gamma(B_R(x_0))}\right)^{\gamma+q-p},
		\end{align*}
		where the positive constant $C$ depends on $n, p, q, s, \alpha, \gamma$ and
$R$.
		\smallskip
		
		\item[(ii)] Let $1<p\le2$ and assume $u$ is locally Lipschitz continuous in $\Omega$. Then
		\begin{align*}
			u\in W^{1+\beta,\gamma}_\mathrm{loc}(\Omega)
		\end{align*}
		for any $\gamma\ge 2$ and $\beta\in \left( 0,\frac{2}{\gamma}\right)$. Moreover, for any $B_{2R}(x_0)\subset\subset\Omega$,
		\begin{align*}
		[\nabla u]^\gamma_{W^{\beta,\gamma}(B_{R/8}(x_0))}\le C\left( 1+\mathfrak{M}_{2R}+\|\nabla u\|_{L^\infty(B_{2R}(x_0))}\right)^{\gamma-2p+2q},
		\end{align*}
where $C$ depends only on $n, s, p, q, \gamma, \beta$ and $R$.
	\end{itemize}
\end{theorem}


The paper is organized as follows. Section \ref{sec2} collects the notation, definitions, and auxiliary results used throughout. Sections \ref{sec3} and \ref{sec4} prove the higher Sobolev regularity of weak solutions to \eqref{main} in the superquadratic range $p\ge2$ and the subquadratic range $1<p\le2$, respectively.
	
\section{Preliminaries}
\label{sec2}

This section introduces the notation and auxiliary results used in the sequel. Throughout the paper, $C$ denotes a positive constant that may change from line to line. Dependencies are displayed when relevant; for example, $C(p,q,n)$ depends only on $p$, $q$, and $n$. We write $B_r(x_0):=\{x\in\mathbb{R}^n:|x-x_0|<r\}$ and abbreviate it to $B_r$ whenever the center is clear or immaterial.

For $1<q<\infty$, the Sobolev space $W^{1,q}(\Omega)$ is defined as
\begin{align*}
	W^{1,q}(\Omega):=\left\lbrace u\in L^q(\Omega): \int_\Omega|\nabla u|^q\,dx<\infty \right\rbrace ,
\end{align*}
and for any $s\in (0,1)$, we also define the fractional Sobolev space $W^{s,q}(\Omega)$ endowed with the norm
\begin{align*}
	\|u\|_{W^{s,q}(\Omega)}:=\|u\|_{L^q(\Omega)}+[u]_{W^{s,q}(\Omega)},
\end{align*}
	where the seminorm $[u]_{W^{s,q}(\Omega)}$ is defined as
\begin{align*}
	[u]_{W^{s,q}(\Omega)}:=\left( \int_\Omega\int_\Omega\frac{|u(x)-u(y)|^q}{|x-y|^{n+s q}}\,dxdy\right)^\frac{1}{q}.
\end{align*}
To ensure that the nonlocal contribution
\begin{align*}
	\int_{\mathbb{R}^n\backslash B_R(x_0)}\frac{|u(y)|^{q-1}}{|x_0-y|^{n+s q}}\,dy
\end{align*}
is finite, we use the tail space
\begin{align*}
	L^{q-1}_{sq}(\mathbb{R}^n):=\left\lbrace u\in L^{q-1}_\mathrm{loc}(\mathbb{R}^n):\int_{\mathbb{R}^n}\frac{|u(x)|^{q-1}}{1+|x|^{n+sq}\,dx<\infty}\right\rbrace.
\end{align*}
For $x_0\in \mathbb{R}^n$ and $R>0$, we define the so-called nonlocal tail
\begin{align*}
	\mathrm{Tail}(u;x_0,R)=\left( R^{sq}\int_{\mathbb{R}^n\backslash B_R(x_0)}\frac{|u|^{q-1}}{|x-x_0|^{n+sq}}\,dx\right)^\frac{1}{q-1}.
\end{align*}
One can easily see that $\mathrm{Tail}(u;x_0,R)$ is well-defined as long as $u\in L^{q-1}_{sq}(\mathbb{R}^n)$.

We also use the monotone function
\begin{align*}
	J_p(t):=|t|^{p-2}t,\quad t\in \mathbb{R}.
\end{align*}
The first- and second-order difference operators are defined, respectively, by
\begin{align*}
	\delta_hu(x):=u_h(x)-u(x):=u(x+h)-u(x)
\end{align*}
and
\begin{align*}
	\delta^2_hu(x):=\delta_h(\delta_hu(x))=u(x+2h)+u(x)-2u(x+h),
\end{align*}
where $h\in\mathbb{R}^n$.

We next define weak solutions. The coefficient $A$ in front of the fractional operator is included to accommodate rescalings. For $a,b\in\mathbb{R}^n$, either $a\cdot b$ or $\langle a,b\rangle$ denotes the Euclidean inner product.
\begin{definition}
	\label{Def3.1}
	Let $p\in (1,\infty)$, $s\in(0,1)$ and $A>0$. We call a function $u\in W^{1,p}_\mathrm{loc}(\Omega)\cap L^{q-1}_{sq}(\mathbb{R}^n)$ a weak solution to
	\begin{align}
		\label{mixa}
		-\Delta_pu+A(-\Delta_q)^su=0\quad\text{in }\Omega,
	\end{align}
	if
	\begin{align*}
		\int_\Omega|\nabla u|^{p-2}\nabla u\cdot\nabla\varphi\,dx+A\int_{\mathbb{R}^n}\int_{\mathbb{R}^n}\frac{J_q(u(x)-u(y))(\varphi(x)-\varphi(y))}{|x-y|^{n+sq}}\,dxdy=0
	\end{align*}
	for every $\varphi\in W^{1,p}(\Omega)$ with compact support in $\Omega$ and extended to $0$ outside $\Omega$.
\end{definition}

We refer to \cite{GL2401} for existence and uniqueness of weak solutions to \eqref{mixa}. The following two Sobolev embeddings, Lemmas \ref{Lem2.19} and \ref{Lem2.20}, give respectively $W^{1,p}\hookrightarrow W^{\gamma,p}$ and $W^{\gamma,p}\hookrightarrow L^{np/(n-\gamma p)}$. For further background on fractional Sobolev spaces, see \cite{DPV12}.

\begin{lemma}
	\label{Lem2.19}
	Let $p\ge 1$ and $\gamma\in (0,1)$. For any $f\in W^{1,p}(B_R)$, we have
	\begin{align*}
		(1-\gamma)\int_{B_R}\int_{B_R}\frac{|f(x)-f(y)|^p}{|x-y|^{n+\gamma p}}\,dxdy\le CR^{(1-\gamma)p}\int_{B_R}|\nabla f|^p\,dx
	\end{align*}
	with the constant $C:=C(n,p)$.
\end{lemma}

\begin{lemma}
	\label{Lem2.20}
	Let $p\ge 1$ and $\gamma\in (0,1)$ such that $\gamma p<n$. For any $f\in W^{\gamma,p}(B_R)$, we have
	\begin{align*}
		\left(\mint_{B_R}|f|^{\frac{np}{n-\gamma p}}\,dx\right)^\frac{n-\gamma p}{np}\le C\left(R^{\gamma p}\int_{B_R}\mint_{B_R}\frac{|f(x)-f(y)|^p}{|x-y|^{n+\gamma p}}\,dxdy+\mint_{B_R}|f|^p\,dx\right)
	\end{align*}
	with the constant $C:=C(n,p,\gamma)$.
\end{lemma}

We state two classical lemmas relating derivatives and finite differences; see, for instance, \cite[Section 7.11]{DT83}.

\begin{lemma}
	\label{Lem3.7}
	Let $q\in (1,\infty)$, $\gamma\in (0,1)$ and $d\in (0,R)$. There exists a constant $C=C(n,q)$ such that for any $w\in W^{1,q}(B_R)$, we have
	\begin{align*}
		\left( \int_{B_{R-d}}|\delta_hw|^q\,dx\right)^{\frac{1}{q}}\le C|h|\|\nabla w\|_{L^q(B_R)} \quad \text{with any } 0<|h|\le d.
	\end{align*}
\end{lemma}

\begin{lemma}
\label{Lem3.8}
	Let $1<p<\infty$, $M>0$ and $0<d<R$. For any $f\in L^p(B_R)$ satisfying
	\begin{align*}
		\int_{B_{R-d}}|\delta_hf|^p\,dx\le M^p|h|^p
	\end{align*}
	for any $0<|h|\le d$, $f$ is weakly differentiable in $B_{R-d}$. Moreover,
	\begin{align*}
		\int_{B_{R-d}}|\nabla f|^p\,dx\le C(n)M^p.
	\end{align*}
\end{lemma}

The next lemma, adapted from \cite[Lemma 2.9]{DKLN24}, will be used to establish the fractional differentiability of gradients.

\begin{lemma}
	\label{Lem2.18}
	For fixed $q\in [1,\infty)$, $R>0$ and $h_0\in (0,R)$. Let $g\in W^{1,q}(B_{R+6h_0}(x_0))$ satisfy
	\begin{align*}
		\left( \sup\limits_{0<|h|<h_0}\int_{B_{R+4h_0}(x_0)}\frac{|\delta_h^2g|^q}{|h|^{q(1+\gamma)}}\,dx\right)^\frac{1}{q}<M
	\end{align*}
	for some constants $M>0$ and $\gamma\in(0,1)$. Then for any $\tilde{\gamma}\in (0,\gamma)$, we have
	\begin{align*}
		[\nabla g]^q_{W^{\tilde{\gamma},q}(B_R(x_0))}&\le \frac{Ch_0^{q(\gamma-\tilde{\gamma})}M^q}{(\gamma-\tilde{\gamma})\gamma^q(1-\gamma)^q}\\
		&\quad +\frac{Ch_0^{q(\gamma-\tilde{\gamma})}}{\tilde{\gamma}(\gamma-\tilde{\gamma})\gamma^q(1-\gamma)^q}\frac{(R+4h_0)^{q+n}}{h_0^{q(1+\gamma)}}\int_{B_{R+4h_0}}|\nabla g|^q\,dx,
	\end{align*}
	where $C$ depends only on $n, q$.
\end{lemma}

\section{Higher regularity in the superquadratic case}
\label{sec3}

In this section, we prove that $|\nabla u|^{\frac{p-2}{2}}\nabla u\in W^{1,2}_\mathrm{loc}(\Omega;\mathbb{R}^n)$ and that $\nabla u$ is locally almost $W^{2/\gamma,\gamma}$. These are precisely Theorems \ref{th1}(i) and \ref{th2}(i).

\subsection{Regularity for nonlinear expressions of gradients}

We begin by differentiating \eqref{main} at the discrete level to obtain a difference-quotient estimate for $|\nabla u|^{\frac{p-2}{2}}\nabla u$. This estimate directly implies its weak differentiability.
\begin{lemma}[Energy estimate]
	\label{Lem2}
	Let $p\ge2$, $s\in (0,1)$ and $u$ be a weak solution of \eqref{main} in $\Omega$. For any $0<r<R$, $d\in (0,\frac{R-r}{2}]$ with $B_{R+d}:=B_{R+d}(x_0)\subset\subset\Omega$, and $0<|h|\le d$, it holds that
	\begin{align*}
		&\quad\int_{B_r}\frac{|F(x+h)-F(x)|^2}{|h|^2}\,dx\\
		&\le \frac{C}{(R-r)^p}\int_{B_R}\frac{|\delta_hu|^p}{|h|^p}dx+C\int_{B_R}|\nabla u|^p\,dx\\
		&\quad+\frac{C\mathfrak{M}^{q-p}_{R}}{(R-r)^2}\left(\int_{B_{R+d}}\int_{B_{R+d}}\frac{|u(x)-u(y)|^p}{|x-y|^{n+sq}}\,dxdy\right)^{\frac{p-2}{p}}\left(\frac{R^{p-sq}}{p-sq}\int_{B_R}\frac{|\delta_hu|^p}{|h|^p}\,dx\right)^\frac{2}{p}\\
		&\quad+\frac{C}{R^{sq}}\left( \frac{R}{R-r}\right)^{n+sq}\mathfrak{M}_R^{q-2}\int_{B_R}\frac{|\delta_hu|^2}{|h|^2}\,dx\\
		&\quad+\frac{C}{R^{sq+1}}\left( \frac{R}{R-r}\right)^{n+sq+1}\mathfrak{M}_R^{q-1}\int_{B_R}\frac{|\delta_hu|}{|h|}\,dx,
	\end{align*}
	where $C\ge1$ depends on $n,p,q,s$ and $F$ is defined as
	\begin{align*}
		F(x)=|\nabla u(x)|^\frac{p-2}{2}\nabla u(x). 
	\end{align*}
\end{lemma}

\begin{proof}
	We first set $\tilde{r}:=r+\frac{d}{2}$ and $\tilde{R}:=R-\frac{d}{2}$. Let $\eta\in C^\infty_0(B_{\frac{R+r}{2}}(x_0))$ be a cut-off function so that $0\le\eta\le1$, $\eta=1$ in $B_{\tilde{r}}$ and $|\nabla \eta|\le\frac{C}{R-r}$. By testing the problem \eqref{main} with $\varphi(x)$ and $\varphi_{-h}(x):=\varphi(x-h)$ respectively, where $\varphi\in W^{1,p}(B_R)$ with $\mathrm{supp}\,\varphi\subset B_{\frac{R+r}{2}}$, we have
	\begin{align*}
		\int_{B_{\tilde{R}}}|\nabla u|^{p-2}\nabla u\cdotp\nabla\varphi\,dx+\int_{\mathbb{R}^n}\int_{\mathbb{R}^n}J_q(u(x)-u(y))(\varphi(x)-\varphi(y))\,d\mu=0
	\end{align*}
	and
	\begin{align*}
		\int_{B_{\tilde{R}}}|\nabla u_h|^{p-2}\nabla u_h\cdotp\nabla\varphi\,dx+\int_{\mathbb{R}^n}\int_{\mathbb{R}^n}J_q(u_h(x)-u_h(y))(\varphi(x)-\varphi(y))\,d\mu=0.
	\end{align*}
Here and in the sequel we denote $d\mu=\frac{dxdy}{|x-y|^{n+sq}}$. Subtracting the two identities above and choosing the test function
	\begin{align*}
		\varphi(x)=\eta^2(x)\frac{u(x+h)-u(x)}{|h|^2},
	\end{align*}
	we have
	\begin{align}
		\label{2.1}
			0
			&=\int_{B_{\tilde{R}}}|h|^{-2}\eta^2\left( |\nabla u_h|^{p-2}\nabla u_h-|\nabla u|^{p-2}\nabla u\right)\cdot(\nabla u_h-\nabla u)\,dx \nonumber\\
			&\quad+2\int_{B_{\tilde{R}}}|h|^{-2}\eta (u_h-u)\left( |\nabla u_h|^{p-2}\nabla u_h-|\nabla u|^{p-2}\nabla u\right)\cdot\nabla \eta \,dx \nonumber\\
			&\quad+\int_{B_{\tilde{R}}}\int_{B_{\tilde{R}}}\left(J_q(u_h(x)-u_h(y))-J_q(u(x)-u(y))\right)(\varphi(x)-\varphi(y))\,d\mu \nonumber\\
			&\quad+2\int_{\mathbb{R}^n\backslash B_{\tilde{R}}}\int_{B_{\frac{R+r}{2}}}\left(J_q(u_h(x)-u_h(y))-J_q(u(x)-u(y))\right)\varphi(x)\,d\mu \nonumber\\
			&=:I_1+2I_2+I_3+2I_4.
	\end{align}

First, the basic inequality \cite[page 99, Inequality (V)]{L19} gives
	\begin{align}
		\label{2.2}
		I_1\ge \int_{B_{\tilde{R}}}\frac{4}{p^2}\frac{\eta^2| F(x+h)-F(x)|^2}{|h|^2}\,dx.
	\end{align}
	Additionally, we utilize Young's inequality to derive
	\begin{align}
		\label{2.3}
			|I_2|&\le C(p)\int_{B_{\tilde{R}}}\frac{|\delta_h u||\nabla \eta|\left( |\nabla u_h|^\frac{p-2}{2}+|\nabla u|^\frac{p-2}{2}\right)}{|h|}\frac{\eta\left| |\nabla u_h|^{\frac{p-2}{2}}\nabla u_h-|\nabla u|^{\frac{p-2}{2}}\nabla u\right| }{|h|}\,dx \nonumber\\
			&\le \varepsilon\int_{B_{\tilde{R}}}\eta^2\frac{\left| |\nabla u_h|^{\frac{p-2}{2}}\nabla u_h-|\nabla u|^{\frac{p-2}{2}}\nabla u\right|^2}{|h|^2}\,dx \nonumber\\
			&\quad+C(\varepsilon,p)\int_{B_{\tilde{R}}}\left| \frac{\delta_hu}{|h|}\right|^p|\nabla \eta|^pdx+C(\varepsilon,p)\int_{B_{\tilde{R}}}\left( |\nabla u_h|^{\frac{p-2}{2}}+|\nabla u|^\frac{p-2}{2}\right)^\frac{2p}{p-2}\chi_{\mathrm{supp}\,\eta}\,dx \nonumber\\
			&\le \varepsilon\int_{B_{\tilde{R}}}\eta^2\frac{\left| |\nabla u_h|^{\frac{p-2}{2}}\nabla u_h-|\nabla u|^{\frac{p-2}{2}}\nabla u\right|^2}{|h|^2}\,dx \nonumber\\
			&\quad+C(\varepsilon,p)\int_{B_{\tilde{R}}}\left|\frac{\delta_hu}{|h|}\right|^p|\nabla \eta|^pdx+C(\varepsilon,p)\int_{B_{\tilde{R}}}|\nabla u|^p\,dx.
	\end{align}
We next estimate $I_3$ and $I_4$. Lemma 2.3 in \cite{BDL2402} gives
	\begin{align}
		\label{2.4}
			I_3
&\ge \frac{1}{C(q)}\int_{B_{\tilde{R}}}\int_{B_{\tilde{R}}}\left( |W_h(x,y)|+|W_0(x,y)|\right)^{q-2}\frac{|\delta_hu(x)-\delta_hu(y)|^2}{|h|^2}(\eta^2(x)+\eta^2(y))\,d\mu \nonumber\\
			&\quad -C(q)\int_{B_{\tilde{R}}}\int_{B_{\tilde{R}}}\left( |W_h(x,y)|+|W_0(x,y)|\right)^{q-2}\left( |\delta_hu(x)|+|\delta_hu(y)|\right) ^2\frac{|\eta(x)-\eta(y)|^2}{|h|^2}\,d\mu \nonumber\\
			&\ge \frac{1}{C(q)}\int_{B_{\tilde{R}}}\int_{B_{\tilde{R}}}|\delta_hu(x)-\delta_hu(y)|^q(\eta^2(x)+\eta^2(y))\,d\mu \nonumber\\
			&\quad -C(q)\int_{B_{\tilde{R}}}\int_{B_{\tilde{R}}}\left( |W_h(x,y)|+|W_0(x,y)|\right)^{q-2} \nonumber \\
			&\qquad\times \left( |\delta_hu(x)|+|\delta_hu(y)|\right) ^2\frac{|\eta(x)-\eta(y)|^2}{|h|^2}\,d\mu,
	\end{align}
	where we abbreviated $W_h(x,y):=u(x+h)-u(y+h)$. On the other hand, by H\"{o}lder's inequality and the properties of $\eta$,
	\begin{align}
		\label{2.5}
			&\quad\int_{B_{\tilde{R}}}\int_{B_{\tilde{R}}}\left( |W_h(x,y)|+|W_0(x,y)|\right)^{q-2}\left( |\delta_hu(x)|+|\delta_hu(y)|\right) ^2\frac{|\eta(x)-\eta(y)|^2}{|h|^2}\,d\mu\nonumber\\
			&\le \frac{C(n)\mathfrak{M}_{2R}^{q-p}}{(R-r)^2}\int_{B_R}\int_{B_R}\frac{\left( |W_h(x,y)|+|W_0(x,y)|\right)^{p-2}\left( |\delta_hu(x)|+|\delta_hu(y)|\right) ^2}{|h|^2|x-y|^{n+sq-2}}\,dxdy\nonumber\\
			&\le \frac{C(n)\mathfrak{M}_{2R}^{q-p}}{(R-r)^2}\left( \int_{B_R}\int_{B_R}\frac{\left( |W_h(x,y)|+|W_0(x,y)|\right)^p}{|x-y|^{n+sq}}\,dxdy\right)^{\frac{p-2}{p}}\nonumber \\
			&\qquad\qquad\qquad\times\left(\int_{B_R}\int_{B_R} \frac{\left( |\delta_hu(x)|+|\delta_hu(y)|\right) ^p}{|h|^p|x-y|^{n+sq-p}}\,dxdy\right)^{\frac{2}{p}}\nonumber \\
			&\le \frac{C(n,p)\mathfrak{M}_{2R}^{q-p}}{(R-r)^2}\left( \int_{B_{R+d}}\int_{B_{R+d}}\frac{|u(x)-u(y)|^p}{|x-y|^{n+sq}}\,dxdy\right)^{\frac{p-2}{p}}\left(\frac{R^{p-sq}}{p-sq}\int_{B_R}\left| \frac{\delta_hu}{|h|}\right|^p\,dx \right)^{\frac{2}{p}}.
	\end{align}
For the nonlocal integral $I_4$, it follows from the finite-difference tail estimate \cite[Lemma 3.1]{BDL2402} that
	\begin{align}
		\label{2.6}
			|I_4|&\le \frac{C}{R^{sq}}\left( \frac{R}{R-r}\right)^{n+sq}\mathfrak{M}_R^{q-2}\int_{B_R}\frac{|\delta_hu|^2}{|h|^2}\,dx \nonumber\\
			&\qquad+\frac{C|h|}{R^{sq+1}}\left( \frac{R}{R-r}\right) ^{n+sq+1}\mathfrak{M}_R^{q-1}\int_{B_R}\frac{|\delta_hu|}{|h|^2}\,dx  \nonumber\\
			&=\frac{C}{R^{sq}}\left( \frac{R}{R-r}\right)^{n+sq}\mathfrak{M}_R^{q-2}\int_{B_R}\frac{|\delta_hu|^2}{|h|^2}\,dx \nonumber\\
			&\qquad+\frac{C}{R^{sq+1}}\left( \frac{R}{R-r}\right)^{n+sq+1}\mathfrak{M}_R^{q-1}\int_{B_R}\frac{|\delta_hu|}{|h|}\,dx.
	\end{align}

	Finally, combining \eqref{2.2}--\eqref{2.6} with \eqref{2.1} and taking $\varepsilon=\frac{2}{p^2}$ yields
	\begin{align*}
		&\quad\int_{B_R}\eta^2\frac{\left| F(x+h)-F(x)\right|^2 }{|h|^2}\,dx\\
		&\le C\int_{B_R}|\nabla \eta |^p\left|\frac{\delta_h u}{|h|}\right|^p\,dx+C\int_{B_R}|\nabla u|^p\,dx\\
		&\quad+\frac{C\mathfrak{M}_{2R}^{q-p}}{(R-r)^2}\left(\int_{B_{R+d}}\int_{B_{R+d}}\frac{|u(x)-u(y)|^p}{|x-y|^{n+sq}}\,dxdy\right)^{\frac{p-2}{p}}\left(\frac{R^{p-sq}}{p-sq}\int_{B_R}\left|\frac{\delta_hu}{|h|}\right|^p\,dx \right)^{\frac{2}{p}}\\
		&\quad+\frac{C}{R^{sq}}\left( \frac{R}{R-r}\right)^{n+sq}\mathfrak{M}_R^{q-2}\int_{B_R}\frac{|\delta_hu|^2}{|h|^2}\,dx+\frac{C}{R^{sq+1}}\left( \frac{R}{R-r}\right) ^{n+sq+1}\mathfrak{M}_R^{q-1}\int_{B_R}\frac{|\delta_hu|}{|h|}\,dx,
	\end{align*}
	where $C\ge1$ depends on $n,p,q,s$.
\end{proof}

Now we are ready to prove Theorem \ref{th1} (i) from the previous energy inequalities.

\smallskip

\begin{proof}[Proof of Theorem \ref{th1} \rm (i)]
	We take $d=\frac{R-r}{2}$. Then it follows from Lemma \ref{Lem2} that
	\begin{align*}
		&\quad\int_{B_r}\frac{|F(x+h)-F(x)|^2}{|h|^2}\,dx\\
		&\le C\left( 1+\frac{1}{(R-r)^p}\right)\int_{B_R}|\nabla u|^p\,dx\\
		&\quad+\frac{C\mathfrak{M}^{q-p}_{2R}}{(R-r)^2}\left[ \frac{R^{p-sq}}{p-sq}\int_{B_{2R}}|\nabla u|^p\,dx\right]^{\frac{p-2}{p}} \left[ \frac{R^{p-sq}}{p-sq}\int_{B_{2R}}|\nabla u|^p\,dx\right]^{\frac{2}{p}}\\
		&\quad+\frac{C}{R^{sq}}\left( \frac{R}{R-r}\right)^{n+sq}\mathfrak{M}_R^{q-2}\int_{B_{2R}}|\nabla u|^2\,dx\\
		&\quad+\frac{C}{R^{sq+1}}\left( \frac{R}{R-r}\right) ^{n+sq+1}\mathfrak{M}_R^{q-1}\int_{B_{2R}}|\nabla u|\,dx
	\end{align*}
	for any $0<|h|\le d=\frac{R-r}{2}$, where we have used Lemmas \ref{Lem2.19} and \ref{Lem3.7}. Since $F=|\nabla u|^{\frac{p-2}{2}}\nabla u$ belongs to $L^2(B_R)$, Lemma \ref{Lem3.8} implies
	\begin{align*}
		\int_{B_r}|\nabla F|^2\,dx&\le C\left( 1+\frac{1}{(R-r)^p}\right)\int_{B_{2R}}|\nabla u|^p\,dx+\frac{CR^{p-sq}\mathfrak{M}_{2R}^{q}}{(R-r)^2}\\
		&\quad+\frac{CR^{p-sq}}{(R-r)^2}\left( \int_{B_{2R}}|\nabla u|^p\,dx\right) ^\frac{q}{p}+\frac{CR^{n\left(2-\frac{2}{p}\right) }}{(R-r)^{n+sq}}\left( \int_{B_{2R}}|\nabla u|^p\,dx\right) ^{\frac{q}{p}}\\
		&\quad+\frac{CR^{n\left(2-\frac{1}{p}\right) }}{(R-r)^{n+sq+1}}\left( \int_{B_{2R}}|\nabla u|^p\,dx\right) ^\frac{q}{p}\\
		&\quad+\left( \frac{CR^{n\left(2-\frac{2}{p}\right) }}{(R-r)^{n+sq}}+\frac{CR^{n\left(2-\frac{1}{p}\right) }}{(R-r)^{n+sq+1}}\right)\mathfrak{M}_R^q.
	\end{align*}
	By taking $r=\frac{R}{2}$, this becomes
	\begin{align*}
		\int_{B_{\frac{R}{2}}}|\nabla F|^2dx\le C(n,p,q,s,R)\left(1+ \|\nabla u\|_{L^p(B_{2R})}+\mathfrak{M}_{2R}\right)^q.
	\end{align*}
The proof is finished.
\end{proof}


\subsection{Fractional Sobolev regularity for gradients} \label{sec3.2}

For clarity, we first reduce the problem by rescaling to the normalized regime $\|u\|_{L^\infty(B_1)}\le1$ and $\mathrm{Tail}(u;0,1)\le1$. We establish higher fractional differentiability in this regime and then scale back to obtain the general estimates.

\subsubsection{Reduction of the problem}
Assume that $1<p<\infty$ and $u$ is a locally bounded weak solution to \eqref{main} in $\Omega$.
Let $B_{2R}(x_0)\subset\subset \Omega$.

Moreover, let
\begin{align*}
	\tilde{u}(x)=\frac{1}{\mathfrak{M}_R}u(Rx+x_0).
\end{align*}
Since the desired estimates are local, it suffices to work on a ball. Thus, if $u$ is a weak solution of \eqref{main} in $B_{2R}(x_0)$, then $\tilde{u}$ is a weak solution of
\begin{align*}
	-\Delta_p\tilde u+R^{p-sq}\mathfrak{M}_R^{q-p}(-\Delta_q)^s\tilde u=0\quad\text{in }B_2(0).
\end{align*}
In addition, we know that
\begin{align*}
	\|\tilde{u}\|_{L^\infty(B_1)}\le 1\quad \text{and}\quad \left( \int_{\mathbb{R}^n\backslash B_{1}}\frac{|\tilde{u}|^{q-1}}{|x|^{n+sq}}\,dx\right)^\frac{1}{q-1}\le 1.
\end{align*}
Because the coefficient $R^{p-sq}\mathfrak{M}_R^{q-p}$ appears in front of $(-\Delta_q)^s$, we consider the mixed local and nonlocal problem
\begin{align}
		\label{main_b}
		-\Delta_pu+A(-\Delta_q)^su=0 \quad\text{in } B_2(0)
	\end{align}
with $A>0$.

\subsubsection{Fractional differentiability}


The next lemma provides an energy estimate for the second-order differences of solutions, from which we derive the fractional differentiability of the gradient.

\begin{lemma}
	\label{Lem1}
	Let $p\ge 2$, $0<s<1$ and $A>0$. Suppose that $u$ is a weak solution of \eqref{main_b}
	with
	\begin{align*}
		\|u\|_{L^\infty(B_1)}\le 1\quad \text{and}\quad \left( \int_{\mathbb{R}^n\backslash B_{1}}\frac{|u|^{q-1}}{|x|^{n+sq}}\,dx\right)^\frac{1}{q-1}\le 1.
	\end{align*}
	For any $0<r<R$, $d\in (0,\frac{R-r}{5}]$ with $B_{R+d}\subset\subset B_1$ and $0<|h|\le d$, we have
	\begin{align*}
		\int_{B_r}|\delta_h^2u|^\gamma\,dx\le C(1+A)\frac{|h|^{\gamma+2}}{(R-r)^\beta}\left( \int_{B_{R+d}}|\nabla u|^\gamma\,dx+1\right)
	\end{align*}
	provided $u\in W^{1,\gamma}(B_1)$ with $\gamma\ge p$. Here $\beta:=\max\left\lbrace q,n+sq+1\right\rbrace$ and $C=C(n,p,q,\gamma,s)$.
\end{lemma}

\begin{proof}
	Take a cut-off function $\eta\in C^\infty_0\left( B_{\frac{R+r}{2}}\right) $ such that
	\begin{align*}
		\eta\equiv1\quad \text{in }B_r,\ \eta=0\quad \text{in }\mathbb{R}^n\backslash B_{\frac{R+r}{2}}\quad\text{and} \quad |\nabla\eta|\le \frac{C}{R-r}\quad \text{in }B_{\frac{R+r}{2}}.
	\end{align*}
	For $\alpha\ge 1$ and $0<|h|\le d$, define
	\begin{align*}
		\phi(x)=\frac{J_{\alpha+1}(u_h(x)-u(x))\eta^q(x)}{|h|^{1+\alpha}}.
	\end{align*}
	Testing \eqref{main_b} by the functions $\phi$ and $\phi_{-h}:=\phi(x-h)$, and arranging the weak formulation carefully, we get
	\begin{align}
		\label{1.1}
			0&=\int_{B_R}\left\langle |\nabla u_h|^{p-2}\nabla u_h-|\nabla u|^{p-2}\nabla u,\nabla \left( J_{\alpha+1}(u_h-u)\eta^q\right) \right\rangle |h|^{-1-\alpha}\,dx  \nonumber\\
			&\quad+A\int_{B_R}\int_{B_R}|h|^{-1-\alpha}\left( J_q(u_h(x)-u_h(y))-J_q(u(x)-u(y))\right)  \nonumber\\
			&\qquad\qquad\qquad\quad\times\left( J_{\alpha+1}(u_h-u)(x)\eta^q(x)-J_{\alpha+1}(u_h-u)(y)\eta^q(y)\right)\,d\mu \nonumber\\
			&\quad+2A\int_{\mathbb{R}^n\backslash B_R}\int_{B_\frac{R+r}{2}}\left( J_q(u_h(x)-u_h(y))-J_q(u(x)-u(y))\right) \nonumber\\
			&\qquad\qquad\qquad\quad\times J_{\alpha+1}(u_h-u)(x)\eta^q(x)|h|^{-1-\alpha}\,d\mu \nonumber\\
			&=:I_1+AI_2+2AI_3.
	\end{align}
From \cite[ Inequality (4.14)]{GL2401}, we find
	\begin{align}
		\label{1.2}
			|h|^{1+\alpha}I_1&\ge C\int_{B_R}\left| \nabla \left( |u_h-u|^\frac{\alpha-1}{p}(u_h-u)\eta^\frac{q}{p}\right) \right|^p\,dx-C\int_{B_R}|u_h-u|^{p+\alpha-1}|\nabla \eta^\frac{q}{p}|^p\,dx\nonumber\\
			&\quad-C\int_{B_R}\left( |\nabla u_h|^{\frac{p-2}{2}}+|\nabla u|^{\frac{p-2}{2}}\right)^2|u_h-u|^{\alpha+1}|\nabla \eta^{\frac{q}{2}}|^2\,dx,
	\end{align}
	where $C$ depends only on $p$ and $\alpha$. By the assumptions $\|u\|_{L^\infty(B_1)}\le 1,\ R+d\le 1$ and Young's inequality,
	\begin{align}
		\label{1.3}
			\int_{B_R}\frac{|\delta_hu|^{p+\alpha-1}}{|h|^{1+\alpha}}|\nabla \eta^\frac{q}{p}|^p\,dx&\le \frac{C}{(R-r)^p}\left( \int_{B_R}|\delta_hu|^\gamma\,dx+\int_{B_R}\left|\frac{\delta_hu}{|h|}\right|^{\frac{\gamma(1+\alpha)}{\gamma-p+2}}\,dx\right)\nonumber\\
			&\le \frac{C}{(R-r)^p}\left( \int_{B_{R+d}}|u|^\gamma\,dx+\int_{B_R}\left|\frac{\delta_hu}{|h|}\right| ^{\frac{\gamma(1+\alpha)}{\gamma-p+2}}\,dx\right)\nonumber\\
			&\le \frac{C}{(R-r)^p}\left( \int_{B_R}\left|\frac{\delta_hu}{|h|}\right|^{\frac{\gamma(1+\alpha)}{\gamma-p+2}}\,dx+1\right)
	\end{align}
	with $C>0$ depending on $n,p,q,\gamma$. Another application of Young's inequality gives
	\begin{align}
		\label{1.4}
			&\quad\int_{B_R}\left( |\nabla u_h|^{\frac{p-2}{2}}+|\nabla u|^{\frac{p-2}{2}}\right)^2  \frac{|\delta_hu|^{1+\alpha}}{|h|^{1+\alpha}}|\nabla \eta^\frac{q}{2}|^2\,dx\nonumber\\
			&\le \frac{C}{(R-r)^2}\left[\int_{B_R}\left( |\nabla u_h|^{\frac{p-2}{2}}+|\nabla u|^{\frac{p-2}{2}}\right)^\frac{2\gamma}{p-2}\,dx+ \int_{B_R}\left|\frac{\delta_hu}{|h|}\right|^{\frac{\gamma(1+\alpha)}{\gamma-p+2}}\,dx\right]\nonumber\\
			&\le \frac{C}{(R-r)^2}\left( \int_{B_{R+d}}|\nabla u|^\gamma\,dx+ \int_{B_R}\left|\frac{\delta_hu}{|h|}\right| ^{\frac{\gamma(1+\alpha)}{\gamma-p+2}}\,dx\right).
	\end{align}
	Here the positive constant $C$ depends on $n,p,q,\gamma$.
	
We next turn to $I_2$. Arguing as in \cite[Step 1, pages 813--817]{BLS18}, we obtain
	\begin{align}
		\label{1.5}
			I_2&\ge C\int_{B_R}\int_{B_R}|h|^{-1-\alpha}\left| |\delta_hu(x)|^{\frac{\alpha-1}{q}}\delta_hu(x)\eta(x)-|\delta_hu(y)|^{\frac{\alpha-1}{q}}\delta_hu(y)\eta(y)\right|^q\,d\mu\nonumber\\
			&\quad-C\int_{B_R}\int_{B_R}|h|^{-1-\alpha}\left(  |\delta_hu(x)|^{q+\alpha-1}+|\delta_hu(y)|^{q+\alpha-1}\right) |\eta(x)-\eta(y)|^q\,d\mu\nonumber\\
			&\quad-C\int_{B_R}\int_{B_R}|h|^{-1-\alpha}\left( |W_h(x,y)|+|W_0(x,y)|\right)^{q-2}\nonumber\\
			&\qquad\qquad\qquad\times\left(  |\delta_hu(x)|+|\delta_hu(y)|\right)^{\alpha+1}|\eta^{\frac{q}{2}}(x)-\eta^{\frac{q}{2}}(y)|^2\,d\mu
	\end{align}
	with $C=C(q,\alpha)$. Using the properties of $\eta$, $\|u\|_{L^\infty(B_1)}\le 1$ and $R\le 1$, we have
	\begin{align}
		\label{1.6}
			&\quad\int_{B_R}\int_{B_R}|h|^{-1-\alpha}\left(  |\delta_hu(x)|^{q+\alpha-1}+|\delta_hu(y)|^{q+\alpha-1}\right) |\eta(x)-\eta(y)|^q\,d\mu\nonumber\\
			&\le \frac{CR^{(1-s)q}}{(q-sq)(R-r)^q}\int_{B_R}\frac{|u_h-u|^{q-p+p+\alpha-1}}{|h|^{1+\alpha}}\,dx\nonumber\\
			&\le \frac{C}{(1-s)(R-r)^q}\left( \int_{B_R}|u_h-u|^\gamma\,dx+\int_{B_R}\left| \frac{\delta_hu}{|h|}\right| ^{\frac{\gamma(1+\alpha)}{\gamma-p+2}}\,dx\right) \nonumber\\
			&\le \frac{C}{(1-s)(R-r)^q}\left( \int_{B_1}|u|^\gamma\,dx+\int_{B_R}\left| \frac{\delta_hu}{|h|}\right| ^{\frac{\gamma(1+\alpha)}{\gamma-p+2}}\,dx\right)\nonumber\\
			&\le \frac{C}{(1-s)(R-r)^q}\left( \int_{B_R}\left| \frac{\delta_hu}{|h|}\right| ^{\frac{\gamma(1+\alpha)}{\gamma-p+2}}\,dx+1\right),
	\end{align}
	where $C=C(n,p,q,\gamma)$ and in the third line we utilized Young's inequality with $(\frac{\gamma}{p-2},\frac{\gamma}{\gamma-p+2})$ for $p>2$. For $p=2$, the case is easier. For the third integral on the right-hand side of \eqref{1.5}, we first observe that
	\begin{align*}
		n+sq-2=\frac{\gamma-p+2}{\gamma}a+\frac{p-2}{\gamma}b
	\end{align*}
	with
	\begin{align*}
		\begin{cases}
			a:=n+sq-p,\\
			b:=n+sq-p+\gamma.
		\end{cases}
	\end{align*}
	Then, applying the properties of $\eta$, Young's inequality, and the embedding in Lemma \ref{Lem2.19}, 
	\begin{align}
		\label{1.7}
			&\quad\int_{B_R}\int_{B_R}\left( |W_h(x,y)|+|W_0(x,y)|\right)^{q-2}\left(  |\delta_hu(x)|+|\delta_hu(y)|\right)^{\alpha+1}\frac{|\eta^{\frac{q}{2}}(x)-\eta^{\frac{q}{2}}(y)|^2}{|h|^{1+\alpha}}\,d\mu \nonumber\\
			&\le \frac{C}{(R-r)^2}\int_{B_R}\int_{B_R}\frac{\left( |W_h(x,y)|+|W_0(x,y)|\right)^{q-p+p-2}}{|x-y|^{\frac{(p-2)b}{\gamma}}}\frac{\left( |\delta_hu(x)|+|\delta_hu(y)|\right)^{\alpha+1}}{|h|^{1+\alpha}|x-y|^{\frac{(\gamma-p+2)a}{\gamma}}}\,dxdy  \nonumber\\
			&\le \frac{C}{(R-r)^2}\Bigg( \int_{B_R}\int_{B_R}\frac{\left( |W_h(x,y)|+|W_0(x,y)|\right)^{\gamma}}{|x-y|^{b}}\,dxdy \nonumber\\
			&\qquad\qquad\qquad+\int_{B_R}\int_{B_R}\frac{1}{|x-y|^{a}}\left( \frac{|\delta_hu(x)|+|\delta_hu(y)|}{|h|}\right)^\frac{\gamma(\alpha+1)}{\gamma-p+2}\,dxdy \Bigg)  \nonumber\\
			&\le \frac{C}{(R-r)^2}\left( \int_{B_{R+d}}\int_{B_{R+d}}\frac{|u(x)-u(y)|^\gamma}{|x-y|^{b}}\,dxdy+\frac{1}{p-sq}\int_{B_R}\left|\frac{\delta_hu}{|h|}\right| ^{\frac{\gamma(1+\alpha)}{\gamma-p+2}}\,dx\right) \nonumber\\
			&\le \frac{C}{(R-r)^2}\left( \frac{1}{p-sq}\int_{B_{R+d}}|\nabla u|^\gamma\,dx+\frac{1}{p-sq}\int_{B_R}\left| \frac{\delta_hu}{|h|}\right| ^{\frac{\gamma(1+\alpha)}{\gamma-p+2}}\,dx\right)  \nonumber\\
			&\le \frac{C}{(1-s)(R-r)^2}\left( \int_{B_{R+d}}|\nabla u|^\gamma\,dx+\int_{B_R}\left| \frac{\delta_hu}{|h|}\right| ^{\frac{\gamma(1+\alpha)}{\gamma-p+2}}\,dx\right),
	\end{align}
	where we exploited $R+d\le 1$ and the positive constant $C$ depends on $n,p,q,\alpha$.
	
We next estimate the nonlocal integral $I_3$. An inspection of the proof of \cite[Lemma 3.1]{BDL2402} gives
	\begin{align*}
		&\quad\left| \int_{\mathbb{R}^n\backslash B_R}\frac{J_q(u_h(x)-u_h(y))-J_q(u(x)-u(y))}{|x-y|^{n+sq}}\,dy\right|\\
		&\le C\frac{\mathfrak{m}^{q-2}}{R^{sq}}\left( \frac{R}{R-r}\right)^{n+sq}|\delta_hu(x)|+C\frac{|h|\mathfrak{m}^{q-1}}{R^{sq+1}} \left( \frac{R}{R-r}\right)^{n+sq+1}\\
		&\quad+C\frac{|h|\|u\|^{q-1}_{L^\infty(B_{R+d})}}{R^{sq+1}}\left( \frac{R}{R-r}\right)^{n+sq}
	\end{align*}
	for any $x\in B_{\frac{R+r}{2}}$. Here $C>0$ depends only on $n,q,s$ and $\mathfrak{m}$ is given as
	\begin{align*}
		\mathfrak{m}=\|u\|_{L^\infty(B_{R})}+\mathrm{Tail}(u;0,R).
	\end{align*}
	By the definition of the tail and the assumptions above,
	\begin{align*}
		\mathrm{Tail}(u;0,R)&\le R^{\frac{sq}{q-1}}\left(\int_{\mathbb{R}^n\backslash B_1}\frac{|u|^{q-1}}{|x|^{n+sq}}\,dx+\int_{B_1\backslash B_R}\frac{|u|^{q-1}}{|x|^{n+sq}}\,dx\right)^{\frac{1}{q-1}}\\
		 &\le R^{\frac{sq}{q-1}}\left[ 1+\|u\|_{L^\infty(B_{1})}\left( \int_{B_1\backslash B_R}|x|^{-n-sq}\,dx\right)^\frac{1}{q-1} \right]\\
		 &\le R^{\frac{sq}{q-1}}\left[1+\left( \frac{\omega_n}{sq}\right)^\frac{1}{q-1}R^{-\frac{sq}{q-1}} \right] \le 1+\left( \frac{\omega_n}{sq}\right)^{\frac{1}{q-1}}.
	\end{align*}
Therefore, $I_3$ can be estimated as
	\begin{align}
		\label{1.8}
			|I_3|&\le \left| \int_{B_{\frac{R+r}{2}}}\int_{\mathbb{R}^n\backslash B_R}\frac{\left( J_q(u_h(x)-u_h(y))-J_q(u(x)-u(y))\right)(|\delta_hu|^{\alpha-1}\delta_hu\eta^q)(x) }{|h|^{1+\alpha}|x-y|^{n+sq}}\,dxdy\right| \nonumber\\
			&\le C\frac{\mathfrak{m}^{q-2}}{R^{sq}}\left( \frac{R}{R-r}\right)^{n+sq}\int_{B_R}\frac{|\delta_hu|^{\alpha+1}}{|h|^{\alpha+1}}\,dx +C\frac{\mathfrak{m}^{q-1}}{R^{sq+1}}\left( \frac{R}{R-r}\right)^{n+sq+1}\int_{B_R}\frac{|\delta_hu|^{\alpha}}{|h|^{\alpha}}\,dx  \nonumber\\
			&\quad+C\frac{\|u\|^{q-1}_{L^\infty(B_{R+d})}}{R^{sq+1}}\left( \frac{R}{R-r}\right)^{n+sq}\int_{B_R}\frac{|\delta_hu|^{\alpha}}{|h|^{\alpha}}\,dx  \nonumber\\
			&\le \frac{C}{R^{sq}}\left( \frac{R}{R-r}\right)^{n+sq+1}\int_{B_R}\frac{|\delta_hu|^{\alpha+1}}{|h|^{\alpha+1}}\,dx+\frac{C}{R^{sq+1}}\left( \frac{R}{R-r}\right)^{n+sq+1}\int_{B_R}\frac{|\delta_hu|^{\alpha}}{|h|^{\alpha}}\,dx   \nonumber\\
			&\le \frac{C}{(R-r)^{n+sq+1}}\left( \int_{B_R}\left| \frac{\delta_hu}{|h|}\right|^{\frac{\gamma(\alpha+1)}{\gamma-p+2}}\,dx+\int_{B_R}\left| \frac{\delta_hu}{|h|}\right|^\frac{\gamma\alpha}{\gamma-p+1}\,dx+1\right) ,
	\end{align}
	where $C$ depends on $n, p, q, s$ and in the last line we have employed Young's inequality and the fact $R\le1$. Combining \eqref{1.2}--\eqref{1.8} with \eqref{1.1} arrives at
	\begin{align*}
		&\quad\int_{B_R}\left| \nabla \left( \frac{|\delta_hu|^{\frac{\alpha-1}{p}}\delta_hu}{|h|^\frac{1+\alpha}{p}}\eta^\frac{q}{p}\right) \right|^p\,dx\\
		&\le \frac{C}{(R-r)^p}\left[\int_{B_R}\left( \frac{|\delta_hu|}{|h|}\right)^\frac{\gamma(1+\alpha)}{\gamma-p+2}\,dx+\int_{B_{R+d}}|\nabla u|^\gamma\,dx+1 \right]\\
		&\quad +A\frac{C}{(p-sq)(R-r)^q} \left[ \int_{B_R}\left( \frac{|\delta_hu|}{|h|}\right)^\frac{\gamma(1+\alpha)}{\gamma-p+2}\,dx+\int_{B_{R+d}}|\nabla u|^\gamma\,dx+1 \right]\\
		&\quad +A\frac{C}{(R-r)^{n+sq+1}}\left[\int_{B_R}\left( \frac{|\delta_hu|}{|h|}\right)^\frac{\gamma(1+\alpha)}{\gamma-p+2}\,dx+\int_{B_{R}}\left( \frac{|\delta_hu|}{|h|}\right)^\frac{\gamma\alpha}{\gamma-p+1}\,dx+1 \right],
	\end{align*}
	where $C=C(n, p, q, \gamma, s)$ and we made use of the facts $2\le p$ and $R-r<1$. Consequently,
	\begin{align*}
		&\quad\int_{B_R}\left| \nabla \left( \frac{|\delta_hu|^{\frac{\alpha-1}{p}}\delta_hu}{|h|^\frac{1+\alpha}{p}}\eta^\frac{q}{p}\right) \right|^p\,dx\\
		&\le \frac{C}{(R-r)^\beta}\left[\int_{B_R}\left( \frac{|\delta_hu|}{|h|}\right)^\frac{\gamma(1+\alpha)}{\gamma-p+2}\,dx+\int_{B_{R+d}}|\nabla u|^\gamma\,dx+1 \right]\\
		&\quad+A\frac{C}{(R-r)^\beta}\left[\int_{B_R}\left( \frac{|\delta_hu|}{|h|}\right)^\frac{\gamma(1+\alpha)}{\gamma-p+2}\,dx+\int_{B_{R}}\left( \frac{|\delta_hu|}{|h|}\right)^\frac{\gamma\alpha}{\gamma-p+1}\,dx+\int_{B_{R+d}}|\nabla u|^\gamma\,dx+1 \right]
	\end{align*}
	with $C=C(n,p,q,\gamma,s)$ and $\beta=\max\left\lbrace q,n+sq+1\right\rbrace $.
	
	Now let $v:=|\delta_hu|^\frac{\alpha-1}{p}\delta_hu\eta^\frac{q}{p}$. Then $v\in W^{1,p}_0(B_R)$. Hence, for any vector $0<|\lambda|\le d$, we obtain
	\begin{align*}
		 \int_{B_{R-d}}|\delta_\lambda v|^p\,dx
&\le |\lambda|^p\int_{B_R}|\nabla v|^p\,dx=|\lambda|^p|h|^{1+\alpha}\int_{B_R}\left| \nabla \left( \frac{v}{|h|^{\frac{1+\alpha}{p}}}\right) \right|^p\,dx\\
		&\le \frac{C|\lambda|^p|h|^{1+\alpha}}{(R-r)^\beta}(1+A)\\
		&\quad \times  \Bigg(\int_{B_R}\left| \frac{\delta_hu}{|h|}\right|^\frac{\gamma(1+\alpha)}{\gamma-p+2}\,dx+\int_{B_{R}}\left|\frac{\delta_hu}{|h|}\right|^\frac{\gamma\alpha}{\gamma-p+1}\,dx+\int_{B_{R+d}}|\nabla u|^\gamma\,dx+1\Bigg).
	\end{align*}
Setting $\lambda=h$ on the left-hand side gives
	\begin{align*}
		\int_{B_{R-d}}|\delta_\lambda v|^p\,dx&\ge \int_{B_r}\left| \delta_\lambda\left( |\delta_hu|^\frac{\alpha-1}{p}\delta_hu\right) \right|^p\,dx \\
		&\xlongequal{\lambda:=h}\int_{B_r}\left| \delta_h\left( |\delta_hu|^\frac{\alpha-1}{p}\delta_hu\right) \right|^p\,dx\\
		&\ge \int_{B_r}|\delta^2_hu|^{p+\alpha-1}\,dx.
	\end{align*}
	We select
	\begin{align*}
		\alpha=\gamma-p+1
	\end{align*}
	and then deduce from Lemma \ref{Lem3.7} that
	\begin{align*}
		\int_{B_r}|\delta^2_hu|^{\gamma}\,dx&\le C\frac{|h|^{\gamma+2}}{(R-r)^\beta}(1+A)\left( \int_{B_{R}}\left( \frac{|\delta_hu|}{|h|}\right)^\gamma\,dx+\int_{B_{R+d}}|\nabla u|^\gamma\,dx+1\right) \\
		&\le C\frac{|h|^{\gamma+2}}{(R-r)^\beta}(1+A)\left( \int_{B_{R+d}}|\nabla u|^\gamma\,dx+1\right).
	\end{align*}
\end{proof}

Lemmas \ref{Lem1} and \ref{Lem2.18} imply that $\nabla u$ belongs to a fractional Sobolev space.

\begin{corollary}
	\label{Cor1}
	Under the assumptions of Lemma \ref{Lem1}, we have $\nabla u\in W^{\alpha,\gamma}(B_\frac{1}{4})$ for every $0<\alpha<\frac{2}{\gamma}$. Moreover,
	\begin{align*}
		[\nabla u]^\gamma_{W^{\alpha,\gamma}(B_\frac{1}{4})}\le C(1+A)\left( \int_{B_1}|\nabla u|^\gamma\,dx+1\right)
	\end{align*}
	holds, where $C>0$ depends only on $n,p,q,\gamma,s$ and $\alpha$.
\end{corollary}

\begin{proof}
	Let $0<r<R\le\frac{1}{2}$. We now apply Lemma \ref{Lem1} with $R:=R$ and $r:=\tilde{r}=\frac{4R+5r}{9}$ to get
	\begin{align*}
		\int_{B_{\widetilde{r}}}|\delta_h^2u|^\gamma\,dx&\le C(1+A)\frac{|h|^{\gamma+2}}{(R-\widetilde{r})^\beta}\left( \int_{B_{R+d}}|\nabla u|^\gamma\,dx+1\right)\\
		&\le C(1+A)\frac{|h|^{\gamma+2}}{(R-r)^\beta}\left( \int_{B_{2R}}|\nabla u|^\gamma\,dx+1\right)
	\end{align*}
	for each $0<|h|\le d=\frac{R-\tilde{r}}{5}=\frac{R-r}{9}$, where $C=C(n,p,q,\gamma,s)$ and $\beta:=\max\left\lbrace q,n+sq+1\right\rbrace $. Observe that $r+4d=\tilde{r}$, $r+6d<R$ and $u\in W^{1,\gamma}(B_R)$. Combining this inequality with Lemma \ref{Lem2.18} implies
	\begin{align*}
		\nabla u\in W^{\alpha,\gamma}(B_r)\quad\text{for any }0<\alpha<\frac{2}{\gamma}.
	\end{align*}

	Furthermore, the estimate below
	\begin{align*}
		&\quad[\nabla u]^\gamma_{W^{\alpha,\gamma}(B_r)}\\
		&\le C(R-r)^{2-\gamma\alpha}\left[ \frac{C(1+A)}{(R-r)^\beta}\left( \int_{B_{2R}}|\nabla u|^\gamma\,dx+1\right)+\frac{1}{\alpha(R-r)^{\gamma+2}}\int_{B_{r+4d}}|\nabla u|^\gamma\,dx\right]\\
		&\le C(1+A)\left[ \left(\frac{1}{(R-r)^{\beta+\gamma\alpha-2}}+\frac{1}{(R-r)^{\gamma+\gamma\alpha}}\right)\int_{B_{2R}}|\nabla u|^\gamma\,dx+\frac{1}{(R-r)^{\beta+\gamma\alpha-2}}\right]\\
		&\le C\frac{1+A}{(R-r)^\theta}\left( \int_{B_{2R}}|\nabla u|^\gamma\,dx+1\right)
	\end{align*}
	holds, where $\theta:=\max\left\lbrace \gamma+\gamma\alpha,\beta+\gamma\alpha-2\right\rbrace $ and $C$ depends on $n,p,q,\gamma,s,\alpha$. In particular, by choosing $r=\frac{1}{2}R$ and $R=\frac{1}{2}$ it follows that
	\begin{align*}
		[\nabla u]^\gamma_{W^{\alpha,\gamma}(B_\frac{1}{4})}\le C(1+A)\left( \int_{B_1}|\nabla u|^\gamma\,dx+1\right)
	\end{align*}
	with $C=C(n, p, q, \gamma, s, \alpha)$.
\end{proof}

We now prove Theorem \ref{th2}(i) by scaling back.

\medskip

\begin{proof}[Proof of Theorem \ref{th2} \rm (i)]
	By the reduction above, the function $\tilde{u}=u(Rx+x_0)/\mathfrak{M}_R$ satisfies the assumptions of Lemma \ref{Lem1}. Hence,
	\begin{align*}
		\left[ \nabla\left( \frac{u(Rx+x_0)}{\mathfrak{M}_R}\right) \right]^\gamma_{W^{\alpha,\gamma}(B_\frac{1}{4})}\le C(1+A)\left( \int_{B_1}\left| \nabla\left( \frac{u(Rx+x_0)}{\mathfrak{M}_R}\right)\right| ^\gamma\,dx+1\right),
	\end{align*}
	namely,
	\begin{align*}
		\left[ \nabla\left( u(Rx+x_0)\right) \right]^\gamma_{W^{\alpha,\gamma}(B_\frac{1}{4})}\le C(1+A)\left( \int_{B_1}\left| \nabla u(Rx+x_0)\right| ^\gamma\,dx+\mathfrak{M}_R^\gamma\right),
	\end{align*}
	where the constant $A=R^{p-sq}\mathfrak{M}_R^{q-p}$ and $B_R(x_0)\subset\subset\Omega$. Now rescaling back, we deduce
	\begin{align*}
		\left[ \nabla u\right]^\gamma_{W^{\alpha,\gamma}(B_\frac{R}{4}(x_0))}&\le C(R^{p-sq-\alpha \gamma}\mathfrak{M}_R^{q-p}+R^{-\alpha \gamma})\left( \int_{B_{R}(x_0)}\left| \nabla u\right| ^\gamma\,dx+R^{n-\gamma}\mathfrak{M}_R^\gamma\right)\\
		&\le C\left(\|\nabla u\|_{L^\gamma(B_R(x_0))}+\mathfrak{M}_R+1\right)^{\gamma+q-p}.
	\end{align*}
	This proves the claim.
\end{proof}

\section{Higher regularity in the subquadratic case}
\label{sec4}

In this section, we first prove higher integrability and then higher differentiability of the gradient for weak solutions of \eqref{main} in the singular range $1<p\le2$. We begin with several auxiliary lemmas from \cite{BDL2401,BDL2402}.
\begin{lemma}
	\label{Lem3.2}
	For any $\gamma>0$ and any $a,b\in \mathbb{R}$, there exist constants $C_1,C_2>0$, depending only on $\gamma$, such that
	\begin{align*}
		C_1\left( |a|+|b|\right)^{\gamma-1}|b-a|\le \left| |a|^{\gamma-1}a-|b|^{\gamma-1}b\right| \le C_2\left( |a|+|b|\right)^{\gamma-1}|b-a|.
	\end{align*}
\end{lemma}

\begin{lemma}
	\label{Lem3.4}
	Let $p\in (1,2]$ and $\delta\ge1$. For any $a,b,c,d\in\mathbb{R}$ and any $e,f>0$, we have
	\begin{align*}
		&\quad\left( J_p(a-b)-J_p(c-d)\right) \left( J_{\delta+1}(a-c)e^p-J_{\delta+1}(b-d)f^p\right)\\
		&\ge \frac{p-1}{2^{\delta+1}}\left( |a-b|+|c-d|\right)^{p-2}\left( |a-c|+|b-d|\right)^{\delta-1}|(a-c)-(b-d)|^2(e^p+f^p)\\
		&\quad-\left( \frac{2^{\delta+1}}{p-1}\right)^{p-1}\left( |a-c|+|b-d|\right)^{p+\delta-1}|e-f|^p.
	\end{align*}
\end{lemma}

\begin{lemma}
	\label{lem3.40}
	Let $p\ge2$ and $\delta\ge1$. Then there exists a constant $C=\widetilde{C}(p)2^\delta$, such that whenever $a,b,c,d\in\mathbb{R}$ and $e,f\in \mathbb{R}_{\ge0}$, we have
	\begin{align*}
		&\quad\left(J_{p}(a-b)-J_{p}(c-d)\right)\left(J_{\delta+1}(a-c)e^2-J_{\delta+1}(b-d)f^2\right)\\
		\ge &\frac{1}{C}\left(|a-b|+|c-d|\right)^{p-2}\left(|a-c|+|b-d|\right)^{\delta-1}\left|(a-c)-(b-d)\right|^{2}\left( e^2+f^2\right)\\
		&\quad-C\left(|a-b|+|c-d|\right)^{p-2}\left(|a-c|+|b-d|\right)^{\delta+1}|e-f|^2.
	\end{align*}
\end{lemma}


The following lemmas provide tail estimates for finite differences and play an important role in the energy estimates.

\begin{lemma}
	\label{Lem3.5}
	Let $q\in (1,2]$ and $s\in (0,1)$. There exists a constant $C=C(n,s,q)$ such that for any $u\in L^{q-1}_{sq}(\mathbb{R}^n)$, $R>0$, $r\in (0,R)$ and $d\in \left( 0,\frac{1}{4}(R-r)\right] $,
	\begin{align*}
		&\quad\left| \int_{\mathbb{R}^n\backslash B_R(x_0)}\frac{J_q(u_h(x)-u_h(y))-J_q(u(x)-u(y))}{|x-y|^{n+sq}}\,dy\right|\\
		&\le \frac{C}{R^{sq}}\left( \frac{R}{R-r}\right)^{n+sq+1}\left( |\delta_hu(x)|^{q-1}+\frac{|h|}{R}\mathfrak{M}_{R+d}^{q-1}\right)
	\end{align*}
	for any $x\in B_{\frac{1}{2}(R+r)}(x_0)$ and any $h\in \mathbb{R}^n\backslash \left\lbrace 0\right\rbrace $ with $0<|h|<d$.
\end{lemma}
\begin{lemma}
	\label{lem3.50}
	Let $q\in [2,+\infty)$ and $s\in(0,1)$. There exists a constant $C=C(n,q,s)$ such that whenever $u\in L^{q-1}_{sq}(\mathbb{R}^n)$, $R>0$, $r\in(0,R)$, and $d\in\left( 0,\frac{1}{4}(R-r)\right]$, we have for any $x\in B_{\frac{1}{2}(R+r)}(x_0)$ and any $h\in\mathbb{R}^n\backslash\left\lbrace 0\right\rbrace$ with $0<|h|<d$ that
	\begin{align*}
		&\quad\left| \int_{\mathbb{R}^n\backslash B_R(x_0)}\frac{J_q\left( u_h(x)-u_h(y)\right)-J_q\left( u(x)-u(y)\right)}{|x-y|^{n+sq}}\,dy\right|\\
		&\le C\frac{|\delta_hu(x)|}{R^{sq}}\left( \frac{R}{R-r}\right)^{n+sq}\mathfrak{M}_{R+d}^{q-2}+C\frac{|h|}{R^{sq+1}}\left(\frac{R}{R-r}\right)^{n+sq+1}\mathfrak{M}_{R+d}^{q-1}.
	\end{align*}
\end{lemma}
\begin{remark}
When we prove the higher differentiability of weak solutions, Lemma \ref{Lem3.5} will be improved under suitable assumptions in the proof of Theorem \ref{th1} (ii).
\end{remark}

\subsection{Higher integrability for the gradients}

As in Subsection \ref{sec3.2}, we continue to consider \eqref{main_b} and study higher integrability of the gradient. We begin with an energy estimate for second-order difference quotients.

\begin{lemma}
	\label{Lem3.6}
	Let $p\in (1,2]$, $s\in (0,1)$ and $A>0$. Suppose that $u\in W^{1,\gamma}(B_1)$ is a weak solution of \eqref{main_b} in $B_2(0)$
	with
	\begin{align}
		\label{3.60}
		\|u\|_{L^\infty(B_1)}\le 1\quad \text{and}\quad \left( \int_{\mathbb{R}^n\backslash B_{1}}\frac{|u|^{q-1}}{|x|^{n+sq}}\,dx\right)^\frac{1}{q-1}\le 1.
	\end{align}
	For any $0<r<R$, $d\in (0,\frac{R-r}{5}]$ with $B_{R+d}\subset B_1$ and $h\in \mathbb{R}^n\backslash\left\lbrace 0\right\rbrace $ with $0<|h|\le d$, we have
	\begin{align*}
		\int_{B_r}|\delta_h^2u|^\gamma\,dx\le C(1+A)^{\frac{p}{2}}\frac{|h|^{\gamma+\frac{p}{2}}}{(R-r)^{(n+sq+1)\frac{p}{2}}}\left( \int_{B_{R+d}}|\nabla u|^\gamma\,dx+1\right),
	\end{align*}
	where $\gamma\ge p$ and the constant $C$ depends on $n, s, p, q,\gamma$.
\end{lemma}

\begin{proof}
	By testing Eq. \eqref{main_b} with $\varphi(x)$ and $\varphi_{-h}(x):=\varphi(x-h)$ respectively, where $\varphi\in W^{1,p}(B_R)$ and $\mathrm{supp}\, \varphi\subset B_{\frac{R+r}{2}}$, we have
	\begin{align}
		\label{3.2}
		\int_{B_2}|\nabla u|^{p-2}\nabla u\cdotp\nabla\varphi\,dx+A\int_{\mathbb{R}^n}\int_{\mathbb{R}^n}\frac{J_q(u(x)-u(y))(\varphi(x)-\varphi(y))}{|x-y|^{n+sq}}\,dxdy=0
	\end{align}
	and
	\begin{align}
		\label{3.4}
		\int_{B_2}|\nabla u_h|^{p-2}\nabla u_h\cdotp\nabla\varphi\,dx+A\int_{\mathbb{R}^n}\int_{\mathbb{R}^n}\frac{J_q(u_h(x)-u_h(y))(\varphi(x)-\varphi(y))}{|x-y|^{n+sq}}\,dxdy=0.
	\end{align}
	Subtracting \eqref{3.2} from \eqref{3.4} and choosing $\varphi=|\delta_hu|^{\gamma-p}\delta_hu\eta^p$, where $\eta\in C^\infty_0\left( B_{\frac{R+r}{2}};[0,1]\right)$, $\eta\equiv1$ in $B_{r+d}$, and $|\nabla \eta|\le \frac{C}{R-r}$ in $B_{\frac{R+r}{2}}$, yields
	\begin{align}
		\label{3.5}
		0
		&= A\int_{\mathbb{R}^n}\int_{\mathbb{R}^n}\left[ J_q(u_h(x)-u_h(y))-J_q(u(x)-u(y))\right] \nonumber\\
		&\qquad\qquad\qquad\times\left( \left[ \eta^p|\delta_hu|^{\gamma-p}\delta_hu\right] (x)-\left[ \eta^p|\delta_hu|^{\gamma-p}\delta_hu\right] (y)\right) \,d\mu\nonumber\\
		&\quad+C\int_{B_R}\left( |\nabla u_h|^{p-2}\nabla u_h-|\nabla u|^{p-2}\nabla u\right) \cdotp \nabla\varphi \,dx=:I_1+I_2
	\end{align}
	with $C=C(p,\gamma)$.
	
	\textbf{Case 1: $q\le2$.} For the fractional part $I_1$, we have
	\begin{align}
		\label{3.6}
		I_1&=A\int_{B_R}\int_{B_R}\left[ J_q(u_h(x)-u_h(y))-J_q(u(x)-u(y))\right]\nonumber\\
		&\qquad\qquad\qquad\times\left( \left[ \eta^p|\delta_hu|^{\gamma-p}\delta_hu\right] (x)-\left[ \eta^p|\delta_hu|^{\gamma-p}\delta_hu\right] (y)\right) \,d\mu\nonumber\\
		&\quad+2A\int_{B_\frac{R+r}{2}}\int_{\mathbb{R}^n\backslash B_R}\left[ J_q(u_h(x)-u_h(y))-J_q(u(x)-u(y))\right]\left[ \eta^p|\delta_hu|^{\gamma-p}\delta_hu\right] (x)\,d\mu\nonumber\\
		&=:I_{1,1}+2I_{1,2}.
	\end{align}
	From Lemma \ref{Lem3.4} with $\delta=\gamma-p+1$, by setting $a=u_h(x)$, $b=u_h(y)$, $c=u(x)$, $d=u(y)$, $e=\eta^{\frac{p}{q}}(x)$ and $f=\eta^{\frac{p}{q}}(y)$, we have
	\begin{align}
		\label{3.7}
		I_{1,1}\ge \frac{1}{C}I_{1,1}^1-CI_{1,1}^2,
	\end{align}
	where $I_{1,1}^1$ and $I_{1,1}^2$ are defined by
	\begin{align*}
		I_{1,1}^1&:=A\int_{B_R}\int_{B_R}\left( |u_h(x)-u_h(y)|+|u(x)-u(y)|\right)^{q-2}\\
		&\qquad\qquad\qquad\times\left( |\delta_hu(x)|+|\delta_hu(y)|\right)^{\gamma-p}|\delta_hu(x)-\delta_hu(y)|^2\left( \eta^p(x)+\eta^p(y)\right)\,d\mu\\
	\end{align*}
	and
	\begin{align*}
		I_{1,1}^2:=A\int_{B_R}\int_{B_R}\left( |\delta_hu(x)|+|\delta_hu(y)|\right)^{q+\gamma-p}|\eta^{\frac{p}{q}}(x)-\eta^{\frac{p}{q}}(y)|^q\,d\mu.
	\end{align*}
	For $I_2$, we have
	\begin{align}
		\label{3.8}
		I_2&=\int_{B_R}\left( |\nabla u_h|^{p-2}\nabla u_h-|\nabla u|^{p-2}\nabla u\right) \cdotp|\delta_hu|^{\gamma-p}\eta^p \nabla \delta_h u\,dx\nonumber\\
		&\quad+\int_{B_R}\left( |\nabla u_h|^{p-2}\nabla u_h-|\nabla u|^{p-2}\nabla u\right) \cdotp \eta^{p-1}|\delta_hu|^{\gamma-p}\delta_hu\nabla\eta  \,dx\nonumber\\
		&=:I_{2,1}+I_{2,2}.
	\end{align}
	By the nonnegativity of $I_{1,1}^1$ and combining \eqref{3.5}--\eqref{3.8}, we find
	\begin{align}
		\label{3.9}
		I_{2,1}\le C\left( |I_{1,1}^2|+|I_{1,2}|\right) +|I_{2,2}|.
	\end{align}
	
	We will estimate all the terms on the right-hand side of \eqref{3.9} one by one.
	First, for $I_{1,1}^2$, since $p\le q$, we have
	\begin{align*}
		|\eta^{\frac{p}{q}}(x)-\eta^{\frac{p}{q}}(y)|^q\le  \frac{C}{(R-r)^q}|x-y|^q.
	\end{align*}
	Moreover, noting $R\le1$ and the boundedness of solutions, we obtain
	\begin{align}
		\label{3.10}
		|I_{1,1}^2|&\le \frac{CA}{(R-r)^q}\int_{B_R}\int_{B_R}\frac{\left( |\delta_hu(x)|+|\delta_hu(y)|\right)^\gamma}{|x-y|^{n+sq-q}}\, dxdy \nonumber \\
		&\le \frac{CA}{(1-s)(R-r)^q}\int_{B_R}|\delta_hu|^\gamma\,dx.
	\end{align}
For $I_{1,2}$, by Lemma \ref{Lem3.5} and the assumption \eqref{3.60},
	\begin{align*}
		|I_{1,2}|&\le \int_{B_\frac{R+r}{2}}\frac{CA}{sR^{sq}}\left( \frac{R}{R-r}\right)^{n+sq+1}\left( |\delta_hu|^{q-1}+\frac{|h|}{R}\right)|\delta_hu|^{\gamma-p+1}\eta^p(x)\,dx\\
		&\le\frac{CA}{s(R-r)^{n+sq+1}}\left(|h|\int_{B_R}|\delta_hu|^{\gamma-p+1}\,dx+\int_{B_R}|\delta_hu|^\gamma\,dx\right).
	\end{align*}
Since $|h|<R<1$ and $p-1\le1$, Young's inequality with exponents $\frac{\gamma}{p-1}$ and $\frac{\gamma}{\gamma-p+1}$ gives
	\begin{align*}
		|h|\int_{B_R}|\delta_hu|^{\gamma-p+1}\,dx&\le C\left(\int_{B_R}|\delta_hu|^\gamma\,dx+|h|^\frac{\gamma}{p-1}\right)\le C\left( \int_{B_R}|\delta_hu|^\gamma\,dx+|h|^\gamma\right).
	\end{align*}
	Hence, we get
	\begin{align}
		\label{3.11}
		|I_{1,2}|\le \frac{CA}{s(R-r)^{n+sq+1}}\left(\int_{B_R}|\delta_hu|^\gamma\,dx+|h|^\gamma\right).
	\end{align}
	For $I_{2,2}$, from H\"{o}lder's inequality with the exponents $\frac{\gamma}{p-1}$ and $\frac{\gamma}{\gamma-p+1}$, we obtain
	\begin{align}
		\label{3.12}
		|I_{2,2}|\le \frac{C}{R-r}\left( \int_{B_\frac{R+r}{2}}|\nabla u|^\gamma\,dx\right) ^\frac{p-1}{\gamma}\left( \int_{B_R}|\delta_hu|^\gamma\,dx\right)^\frac{\gamma-p+1}{\gamma}.
	\end{align}
	
	Next, we will estimate the left-hand side of \eqref{3.9} from below. By using H\"{o}lder's inequality and the following basic inequality
	\begin{align*}
		C(p)\left( |\nabla u_h|+|\nabla u|\right)^{p-2} |\nabla\delta_hu|^2\le \left\langle |\nabla u_h|^{p-2}\nabla u_h-|\nabla u|^{p-2}\nabla u,\nabla u_h-\nabla u\right\rangle,
	\end{align*}
	we have
	\begin{align*}
		&\quad\int_{B_R}|\nabla u_h-\nabla u|^p|u_h-u|^{\gamma-p}\eta^p\,dx\nonumber\\
		&\le \left[ \int_{B_R}\left( |\nabla u_h|+|\nabla u|\right)^{p-2}|\nabla u_h-\nabla u|^2|u_h-u|^{\gamma-p}\eta^p\,dx\right]^\frac{p}{2}\\
		&\qquad\times \left[ \int_{B_R}\left( |\nabla u_h|^p+|\nabla u|^p\right)|u_h-u|^{\gamma-p}\eta^p\, dx\right]^{1-\frac{p}{2}} \\
		&\le C(p)\left[ \int_{B_R}\left\langle |\nabla u_h|^{p-2}\nabla u_h-|\nabla u|^{p-2}\nabla u,\nabla u_h-\nabla u\right\rangle|u_h-u|^{\gamma-p}\eta^p\,dx\right] ^\frac{p}{2}\\
		&\qquad\times \left[ \int_{B_R}\left( |\nabla u_h|^p+|\nabla u|^p\right)|u_h-u|^{\gamma-p}\eta^p\, dx\right]^{1-\frac{p}{2}}\\
		&\le C(p,\gamma)|I_{2,1}|^{\frac{p}{2}}\left( \int_{B_R}|\nabla u|^\gamma\eta^p\,dx\right)^{\frac{p}{\gamma}\left( 1-\frac{p}{2}\right) }\left( \int_{B_R}|\delta_hu|^\gamma\eta^p\,dx\right)^{\frac{\gamma-p}{\gamma}\left( 1-\frac{p}{2}\right) }.
	\end{align*}
	Note that $|\nabla u_h-\nabla u|^p|u_h-u|^{\gamma-p}=\left(\frac{\gamma}{p}\right)^p\left|\nabla \left( |\delta_hu|^{\frac{\gamma}{p}-1} \delta_hu\right) \right|^p$.
	Thus, we infer
	\begin{align}
		\label{3.13}
		&\quad\int_{B_R}\left|\nabla \left( |\delta_hu|^{\frac{\gamma}{p}-1} \delta_hu\right) \eta\right|^p\, dx\nonumber\\
		&\le C(p,\gamma)|I_{2,1}|^{\frac{p}{2}}\left(\int_{B_R}|\nabla u|^\gamma\eta^p\,dx\right)^{\frac{p}{\gamma}\left( 1-\frac{p}{2}\right) }\left(\int_{B_R}|\delta_hu|^\gamma\eta^p\,dx\right)^{\frac{\gamma-p}{\gamma}\left( 1-\frac{p}{2}\right) }\nonumber\\
		&\le C(p,\gamma)\left[ |I_{1,1}^2|^\frac{p}{2}+|I_{1,2}|^\frac{p}{2}+|I_{2,2}|^\frac{p}{2} \right] \left( \int_{B_R}|\nabla u|^\gamma\eta^p\,dx\right)^{\frac{p}{\gamma}\frac{2-p}{2} }  \left( \int_{B_R}|\delta_hu|^\gamma\eta^p\,dx\right)^{\frac{\gamma-p}{\gamma} \frac{2-p}{2} }.
	\end{align}
	Substituting \eqref{3.10}--\eqref{3.12} into \eqref{3.13}, we obtain
	\begin{align*}
		&\quad\int_{B_R}\left|\nabla \left( |\delta_hu|^{\frac{\gamma}{p}-1} \delta_hu\right)\eta \right|^p\, dx\\
		&\le C\left[ \left( \frac{A}{(R-r)^q}\int_{B_R}|\delta_hu|^\gamma\,dx\right)^\frac{p}{2}+\left( \frac{A}{(R-r)^{n+sq+1}}\left[ \int_{B_R}|\delta_hu|^\gamma\,dx+|h|^\gamma\right]\right)^\frac{p}{2}\right.\\
		&\qquad\left.+\left( \frac{1}{R-r}\left( \int_{B_\frac{R+r}{2}}|\nabla u|^\gamma\,dx\right) ^\frac{p-1}{\gamma}\left( \int_{B_R}|\delta_hu|^\gamma\,dx\right)^\frac{\gamma-p+1}{\gamma}\right)^\frac{p}{2}   \right]\\
		&\qquad\times \left(\int_{B_R}|\nabla u|^\gamma\eta^p\,dx\right)^{\frac{p}{\gamma}\left( 1-\frac{p}{2}\right) }\left(\int_{B_R}|\delta_hu|^\gamma\eta^p\,dx\right)^{\frac{\gamma-p}{\gamma}\left( 1-\frac{p}{2}\right) }\\
		&\le C\left( \frac{1+A}{(R-r)^{n+sq+1}}\right)^\frac{p}{2}\left(\int_{B_R}|\nabla u|^\gamma\eta^p\,dx\right)^{\frac{p}{\gamma}\left( 1-\frac{p}{2}\right) }\left(\int_{B_R}|\delta_hu|^\gamma\eta^p\,dx\right)^{\frac{\gamma-p}{\gamma}\left( 1-\frac{p}{2}\right) }\\
		&\qquad\times\left[ \left( \int_{B_R}|\delta_hu|^\gamma\,dx\right)^\frac{p}{2}+|h|^\frac{\gamma p}{2}+\left( \|\nabla u\|_{L^\gamma(B_R)}^{p-1}\left( \int_{B_R}|\delta_hu|^\gamma\,dx\right)^\frac{\gamma-p+1}{\gamma}\right) ^\frac{p}{2} \right],
	\end{align*}
	where $C$ depends only on $n,p,q,s,\gamma$. Since $ |\delta_hu|^{\frac{\gamma}{p}-1} \delta_hu\in W^{1,p}(B_R)$, repeated applications of Lemma \ref{Lem3.7} give, for any $0<|\lambda|\le d$,
	\begin{align*}
		\int_{B_{r}}\left|\delta_\lambda \left( |\delta_hu|^{\frac{\gamma}{p}-1} \delta_hu\right) \right| ^p\,dx&\le |\lambda|^p\int_{B_R}\left|\nabla \left( |\delta_hu|^{\frac{\gamma}{p}-1} \delta_hu\right)\eta \right|^p\, dx\\
		&\le C|\lambda|^p\left( \frac{1+A}{(R-r)^{n+sq+1}}\right)^\frac{p}{2}|h|^{(\gamma-p)\left( 1-\frac{p}{2}\right) } \\
		&\quad\times\left( |h|^\frac{\gamma p}{2}+|h|^{(\gamma-p+1)\frac{p}{2}}\right)\left( \int_{B_{R+d}}|\nabla u|^\gamma\,dx+1\right) \\
		&\le C|\lambda|^p\left( \frac{1+A}{(R-r)^{n+sq+1}}\right)^\frac{p}{2}\left( \int_{B_{R+d}}|\nabla u|^\gamma\,dx+1\right)|h|^{\gamma-\frac{p}{2}}.
	\end{align*}
	Here we have used the fact $|h|<1$. For the left-hand side, we have
	\begin{align*}
		\int_{B_{r}}\left| \delta_\lambda \left( |\delta_hu|^{\frac{\gamma}{p}-1} \delta_hu\right) \right|^p\,dx \xlongequal{\lambda:=h}\int_{B_{r}}\left| \delta_h \left( |\delta_hu|^{\frac{\gamma}{p}-1} \delta_hu\right) \right| ^p\,dx\ge \int_{B_{r}}|\delta^2_hu|^{\gamma}\, dx,
	\end{align*}
	which implies
	\begin{align*}
		\int_{B_r}|\delta_h^2u|^\gamma\,dx\le C(1+A)^{\frac{p}{2}}\frac{|h|^{\gamma+\frac{p}{2}}}{(R-r)^{(n+sq+1)\frac{p}{2}}}\left( \int_{B_{R+d}}|\nabla u|^\gamma\,dx+1\right)
	\end{align*}
	with $C=C(n,s,p,q,\gamma)$.
	
	\textbf{Case 2: }$q>2$. We restrict our attention to the fractional part $I_1$. By Lemma \ref{lem3.40},
	\begin{align}
		\label{3.14}
		I_{1,1}\ge \frac{1}{C}I^3_{1,1}-CI^4_{1,1}
	\end{align}
	with
	\begin{align*}
		I^3_{1,1}&:=A\int_{B_R}\int_{B_R}\left(|u_h(x)-u_h(y)|+|u(x)-u(y)|\right)^{q-2}\left(|\delta_hu(x)|+|\delta_hu(y)|\right)^{\gamma-p}\\
		&\qquad\qquad\qquad\times\left|(\delta_hu(x))-(\delta_hu(y))\right|^{2}\left( \eta^p(x)+\eta^p(y)\right)\,d\mu,
	\end{align*}
	and
	\begin{align*}
		I^4_{1,1}&:=A\int_{B_R}\int_{B_R}\left(|u_h(x)-u_h(y)|+|u(x)-u(y)|\right)^{q-2}\\
		&\qquad\qquad\qquad\times\left( |\delta_hu(x)|+|\delta_hu(y)|\right)^{\gamma-p+2}|\eta^{\frac{p}{2}}(x)-\eta^{\frac{p}{2}}(y)|^2\,d\mu.
	\end{align*}
	Due to $p>sq$ and $|\eta^{\frac{p}{2}}(x)-\eta^{\frac{p}{2}}(y)|^2\le|\eta(x)-\eta(y)|^p\le\frac{C}{(R-r)^p}|x-y|^p$,
	we use the boundedness of $u$ to have
	\begin{align}
		\label{3.15}
		I^4_{1,1}&\le \frac{CA}{(R-r)^p}\int_{B_R}\int_{B_R}\frac{\left(|u_h(x)-u_h(y)|+|u(x)-u(y)|\right)^{q-2}}{|x-y|^{n+sq-p}}\nonumber\\
		&\qquad\quad\qquad\qquad\qquad\times\left( |\delta_hu(x)|+|\delta_hu(y)|\right)^{\gamma-p+2}\,dxdy \nonumber\\
		&\le \frac{C(p,q,\gamma)A}{(p-sq)(R-r)^p}\int_{B_R}|\delta_hu|^{\gamma}\,dx.
	\end{align}
	
	For the nonlocal term $I_{1,2}$, Lemma \ref{lem3.50} gives
	\begin{align*}
		|I_{1,2}|&\le \frac{CA}{R^{sq}}\left( \frac{R}{R-r}\right)^{n+sq+1}\int_{B_{\frac{R+r}{2}}}\left(|\delta_hu|^{\gamma-p+2}+|h||\delta_hu|^{\gamma-p+1}\right)\,dx\\
		&\le \frac{CA}{R^{sq}}\left( \frac{R}{R-r}\right)^{n+sq+1}\left( \int_{B_{R}}|\delta_hu|^\gamma\,dx+|h|^{\frac{\gamma}{p-1}}\right)\\
		&\le \frac{CA}{R^{sq}}\left( \frac{R}{R-r}\right)^{n+sq+1}\left( \int_{B_{R}}|\delta_hu|^\gamma\,dx+|h|^{\gamma}\right).
	\end{align*}
	Hence, by replacing $I^2_{1,1}$ with $I^4_{1,1}$ in \eqref{3.13} and repeating the same computations as in the previous case, we obtain the desired inequality.
\end{proof}

Next, one can easily deduce the following corollaries stating higher differentiability of the gradient by combining Lemma \ref{Lem3.6} and Lemma \ref{Lem2.18}.

\begin{corollary}
	\label{Cor3.11}
	Under the assumptions of Lemma \ref{Lem3.6}, we have $\nabla u\in W^{\alpha,\gamma}(B_\frac{1}{4};\mathbb{R}^n)$ for every $0<\alpha<\frac{p}{2\gamma}$. Moreover, it holds that
	\begin{align*}
		[\nabla u]_{W^{\alpha,\gamma}(B_{\frac{1}{4}})}^\gamma\le C(1+A)^\frac{p}{2}\left( \int_{B_1}|\nabla u|^\gamma\,dx+1\right),
	\end{align*}
	where $C$ depends only on $n,s,p,q,\gamma$ and $\alpha$.
\end{corollary}

\begin{proof}
	Similar to the proof of Corollary \ref{Cor1}, by setting $R=\frac{1}{2}$ and $r=\frac{3}{8}$, one can obtain
	\begin{align*}
		\int_{B_\frac{3}{8}}|\delta_h^2u|^\gamma\,dx\le C(1+A)^{\frac{p}{2}}|h|^{\gamma+\frac{p}{2}}\left( \int_{B_{1}}|\nabla u|^\gamma\,dx+1\right)
	\end{align*}
	with $C=C(n,s,p,q,\gamma)$. Thus, letting $d=\frac{1}{64}$ in Lemma \ref{Lem2.18}, we get $\nabla u\in W^{\alpha,\gamma}(B_\frac{1}{4};\mathbb{R}^n)$ for any $\alpha\in (0,\frac{p}{2\gamma})$. Moreover, the following estimate
	\begin{align*}
		[\nabla u]^\gamma_{W^{\alpha,\gamma}(B_{\frac{1}{4}})}\le C(n,s,p,q,\gamma,\alpha)(1+A)^\frac{p}{2}\left( \int_{B_1}|\nabla u|^\gamma\,dx+1\right)
	\end{align*}
	holds.
\end{proof}
Rescaling the preceding estimate yields a $W^{\alpha,\gamma}$ estimate for the gradient of weak solutions to \eqref{main}.

\begin{corollary}
	\label{Coro3.11}
	Let $p\in (1,2]$ and $s\in(0,1)$. Suppose that $u$ is a weak solution of \eqref{main} satisfying $u\in W^{1,\gamma}_\mathrm{loc}(\Omega)$.
	Then, for any $\alpha\in (0,\frac{p}{2\gamma})$ and $B_R(x_0)\subset\subset\Omega$, we have $u\in W^{1+\alpha,\gamma}\left( B_\frac{R}{4}(x_0)\right)$. Moreover,
	\begin{align*}
		\left[ \nabla u\right]^\gamma_{W^{\alpha,\gamma}(B_\frac{R}{4}(x_0))}\le C\left(1+R^{p-sq}\mathfrak{M}_R^{q-p}\right)^\frac{p}{2}R^{-\alpha \gamma}\left( \int_{B_{R}(x_0)}\left| \nabla u\right| ^\gamma\,dx+R^{n-\gamma}\mathfrak{M}_R^\gamma\right),
	\end{align*}
where $C=C(n,s,p,q,\gamma,\alpha)$.
\end{corollary}

\begin{proof}
	It is easy to see that $\frac{u(Rx+x_0)}{\mathfrak{M}_R}$ fulfills the conditions of Lemma \ref{Lem3.6}. Therefore, from Corollary \ref{Cor3.11} and the reduction in Section \ref{sec3}, we have
	\begin{align*}
		\left[ \nabla\left( \frac{u(Rx+x_0)}{\mathfrak{M}_R}\right) \right]^\gamma_{W^{\alpha,\gamma}(B_\frac{1}{4})}\le C(1+A)^\frac{p}{2}\left( \int_{B_1}\left| \nabla\left( \frac{u(Rx+x_0)}{\mathfrak{M}_R}\right)\right| ^\gamma\,dx+1\right),
	\end{align*}
	that is,
	\begin{align*}
		\left[ \nabla u(Rx+x_0) \right]^\gamma_{W^{\alpha,\gamma}(B_\frac{1}{4})}\le C(1+A)^\frac{p}{2}\left( \int_{B_1}\left| \nabla u(Rx+x_0)\right| ^\gamma\,dx+\mathfrak{M}_R^\gamma\right)
	\end{align*}
	with the constant $A=R^{p-sq}\mathfrak{M}_R^{q-p}$. After rescaling back, we obtain
	\begin{align*}
		\left[ \nabla u\right]^\gamma_{W^{\alpha,\gamma}(B_\frac{R}{4}(x_0))}\le C\left(1+R^{p-sq}\mathfrak{M}_R^{q-p}\right)^\frac{p}{2}R^{-\alpha \gamma}\left( \int_{B_{R}(x_0)}\left| \nabla u\right| ^\gamma\,dx+R^{n-\gamma}\mathfrak{M}_R^\gamma\right).
	\end{align*}
\end{proof}

By the embedding $W^{\alpha,\gamma}\hookrightarrow L^{n\gamma/(n-\alpha\gamma)}$, the fractional differentiability obtained in Corollary \ref{Coro3.11} improves the integrability of the gradient. We use this observation in a Moser-type iteration to prove Theorem \ref{th3.13}.

\smallskip

\begin{proof}[Proof of Theorem \ref{th3.13}]
	We will prove this result by induction. We first set
	\begin{align*}
		\gamma_i:=\left( \frac{n}{n-\frac{p}{4}}\right)^ip,\quad \alpha_i:=\frac{p}{4\gamma_i}\quad\text{for }i\in\mathbb{N}.
	\end{align*}
	Let any $x_0\in \Omega$. There is $R>0$ such that $B_{2R}(x_0)\subset\subset\Omega$. For $i=1$, by applying Corollary \ref{Coro3.11} we obtain
	\begin{align*}
		u\in W^{1+\alpha_1,\gamma_1}\left( B_\frac{R}{4}(x_0)\right).
	\end{align*}
	Thus, due to the boundedness of $u$ and the embedding theorem $W^{\alpha,\gamma}\hookrightarrow L^{\frac{n\gamma}{n-\alpha \gamma}}$, Lemma \ref{Lem2.20}, we arrive at
	\begin{align*}
		u\in W^{1,\gamma_2}\left( B_\frac{R}{4}(x_0)\right),
	\end{align*}
	where $\gamma_2=\frac{n\gamma_1}{n-\alpha_1 \gamma_1}$. Because of the arbitrariness of $x_0$, we conclude that
	\begin{align*}
		u\in W^{1,\gamma_2}_\mathrm{loc}\left( \Omega\right).
	\end{align*}
	
	Next, for any $i\ge2$, by supposing that $u\in W^{1,\gamma_i}_\mathrm{loc}\left( \Omega\right)$ and reduplicating the previous step in this proof, one has
	\begin{align*}
		u\in W^{1,\gamma_{i+1}}_\mathrm{loc}\left( \Omega\right).
	\end{align*}
	At this moment, for any $\gamma\ge p$, noting that
	\begin{align*}
		\lim\limits_{i\rightarrow \infty}\gamma_i=\infty,
	\end{align*}
	there exists $i_0\in \mathbb{N}$ such that $\gamma_{i_0}<\gamma\le \gamma_{i_0+1}$. For any set $E\subset\subset \Omega$, by H\"{o}lder's inequality, we get
	\begin{align*}
		\left( \int_E|\nabla u|^\gamma\,dx\right)^\frac{1}{\gamma}\le  \left( \int_E|\nabla u|^{\gamma_{i_0+1}}\,dx\right)^\frac{1}{\gamma_{i_0+1}}|E|^{\frac{1}{\gamma}-\frac{1}{\gamma_{i_0+1}}}.
	\end{align*}
	That is,
	\begin{align*}
		u\in W^{1,\gamma}_\mathrm{loc}(\Omega)\quad\text{for any }\gamma\in [p,\infty).
	\end{align*}
\end{proof}

\subsection{Higher differentiability for the gradients} \label{sub4.2}


We now prove Theorem \ref{th2}(ii). Under the Lipschitz assumption on $u$, we refine the estimates for the fractional part of the equation beyond Lemma \ref{Lem3.5}. The following integration-by-parts formula for finite differences is essential.

\begin{lemma}
	\label{Lem5.1}
	Let $E\subset\mathbb{R}^k$ be an open set, $F\in L^q(E)$, $\xi\in W^{1,q'}(E)$ with $\mathrm{supp}\,\xi\subset\subset E$ and $h\in
	\mathbb{R}^k$, where $q':=\frac{q}{q-1}$ and $0<|h|<\mathrm{dist}(\mathrm{supp}\,\xi,\partial E)$. Then
	\begin{align*}
		\int_E\xi(x)\delta_hF(x)\,dx=-\int_E\int_0^1F(x+th)\,dt\nabla\xi(x)\cdot h\,dx.
	\end{align*}
\end{lemma}

	In particular, if $\mathbb{R}^k=\mathbb{R}^n\times\mathbb{R}^n$, $E=E_1\times E_2\subset\mathbb{R}^n\times\mathbb{R}^n$, and
$h=(h_1,h_2)\in\mathbb{R}^n\times\mathbb{R}^n$ satisfies $0<|h|<\mathrm{dist}(\mathrm{supp}\,\xi,\partial E)$, then
	\begin{align*}
		\int_{E_1}\int_{E_2}\xi(x,y)\delta_hF(x,y)\,dydx=-\int_{E_1}\int_{E_2}\int_0^1F(x+th_1,y+th_2)\,dt\nabla_{(x,y)}\xi(x,y)\cdot
h\,dydx.
	\end{align*}

\medskip

\begin{proof}[Proof of  Theorem \ref{th2} \rm (ii)]
	For any $B_{2R}:=B_{2R}(x_0)\subset\subset\Omega$ and $\gamma\ge2$, we set
$$
\varphi:=\frac{\eta^2J_\gamma(\delta_hu)}{|h|^\gamma}
$$
	as a test function in the weak formulation to obtain
	\begin{align*}
		\int_\Omega|\nabla u|^{p-2}\nabla u\cdot\nabla\varphi\,dx+\int_{\mathbb{R}^n}\int_{\mathbb{R}^n}\frac{J_q(u(x)-u(y))(\varphi(x)-\varphi(y))}{|x-y|^{n+sq}}\,dxdy=0,
	\end{align*}
	where $h\in \mathbb{R}^n\backslash\left\lbrace 0\right\rbrace $, $0<|h|<\frac{R}{16}$ and the cut-off function $\eta\in C^2_0(B_{\frac{3}{4}R})$ satisfies $\eta\equiv1$ in $B_{\frac{1}{2}R}$, $|\nabla\eta|\le \frac{C(n)}{R}$ and $|D^2\eta|\le \frac{C(n)}{R^2}$.
	
	A direct calculation gives
	\begin{align}
		\label{5.1}
		0&=\int_{B_R}\left( |\nabla u_h|^{p-2}\nabla u_h-|\nabla u|^{p-2}\nabla u\right)\cdot \nabla\varphi \,dx \nonumber\\
		&\quad+\int_{\mathbb{R}^n}\int_{\mathbb{R}^n}\frac{\left(J_q(u_h(x)-u_h(y))-J_q(u(x)-u(y))\right)(\varphi(x)-\varphi(y))}{|x-y|^{n+sq}}\,dxdy \nonumber\\
		&\ge \frac{1}{|h|^\gamma}\int_{B_R}\eta^2|\delta_hu|^{\gamma-2}\left(  |\nabla u_h|^{p-2}\nabla u_h-|\nabla u|^{p-2}\nabla u\right)\cdot(\nabla u_h-\nabla u)\,dx\nonumber\\
		&\quad+2\frac{1}{|h|^\gamma}\int_{B_R}\eta J_\gamma(\delta_hu)\left(  |\nabla u_h|^{p-2}\nabla u_h-|\nabla u|^{p-2}\nabla u\right)\cdot\nabla\eta\,dx\nonumber\\
		&\quad+\int_{B_R}\int_{B_R}\frac{\left(J_q(u_h(x)-u_h(y))-J_q(u(x)-u(y))\right)(\varphi(x)-\varphi(y))}{|x-y|^{n+sq}}\,dxdy \nonumber\\
		&\quad+2\int_{\mathbb{R}^n\backslash B_R}\int_{B_{\frac{7R}{8}}}\frac{\left(J_q(u_h(x)-u_h(y))-J_q(u(x)-u(y))\right)\varphi(x)}{|x-y|^{n+sq}}\,dxdy \nonumber\\
		&=:I_1+2I_2+I_3+2I_4.
	\end{align}
	For $I_3$, thanks to the fact that
	\begin{align*}
		&\quad J_\gamma(\delta_hu(x))\eta^2(x)-J_{\gamma}(\delta_hu(y)) \eta^2(y)\\
		&=\left[ J_\gamma(\delta_hu(x))\eta(x)+J_{\gamma}(\delta_hu(y))\eta(y)\right] (\eta(x)-\eta(y))\\
		&\quad+\left[ J_\gamma(\delta_hu(x))-J_{\gamma}(\delta_hu(y))\right]\eta(x)\eta(y),
	\end{align*}
	we have
	\begin{align}
		\label{5.2}
		I_3&\ge\frac{1}{|h|^\gamma}\int_{B_R}\int_{B_R}\frac{\left(J_q(u_h(x)-u_h(y))-J_q(u(x)-u(y))\right)}{|x-y|^{n+sq}}\nonumber\\
		&\qquad\qquad\qquad\quad\times\left[ J_\gamma(\delta_hu(x))-J_{\gamma}(\delta_hu(y))\right]\eta(x)\eta(y)\,dxdy\nonumber\\
		&\quad+\frac{1}{|h|^\gamma}\int_{B_R}\int_{B_R}\frac{\left(J_q(u_h(x)-u_h(y))-J_q(u(x)-u(y))\right)}{|x-y|^{n+sq}}\nonumber\\
		&\qquad\qquad\qquad\qquad\times\left[ J_\gamma(\delta_hu(x))\eta(x)+J_{\gamma}(\delta_hu(y))\eta(y)\right] (\eta(x)-\eta(y))\,dxdy\nonumber\\
		&=:I_{3,1}+I_{3,2}.
	\end{align}
Since $I_{3,1}\ge0$, combining \eqref{5.1} and \eqref{5.2} gives
	\begin{align}
		\label{5.3}
		I_1\le 2|I_2|+|I_{3,2}|+|I_4|.
	\end{align}
	
	Next, we will estimate all terms on the right-hand side of \eqref{5.3} one by one.

	\smallskip
	
\noindent \textbf{Estimate of $I_{2}$. }Since $u$ is locally Lipschitz continuous and $\mathrm{supp}\,\varphi\subset B_{\frac{3}{4}R}$, integration by parts and Young's inequality give
	\begin{align}
		\label{5.4}
		I_2&=\left| \frac{1}{|h|^\gamma}\int_{B_R}\eta J_\gamma(\delta_hu)\left\langle \left( |\nabla u|^{p-2}\nabla u\right)(x+h)-\left( |\nabla u|^{p-2}\nabla u\right) (x) ,\nabla\eta\right\rangle \,dx\right| \nonumber\\
		&\le \frac{1}{|h|^\gamma}\left| \int_{B_R}\left\langle \int_0^1\left( |\nabla u|^{p-2}\nabla u\right)(x+th)\,dt,\left[ \nabla \left( \eta J_\gamma(\delta_hu)\nabla \eta\right) \right]\cdot h  \right\rangle \,dx\right| \nonumber\\
		&=\frac{1}{|h|^\gamma}\left| \int_{B_R}\left\langle \int_0^1\left( |\nabla u|^{p-2}\nabla u\right)(x+th)\,dt,\left[ J_\gamma(\delta_hu)\left( \nabla\eta\otimes\nabla\eta+\eta D^2\eta\right) \right] \cdot h  \right\rangle \,dx\right| \nonumber\\
		&\quad+\frac{\gamma-1}{|h|^\gamma}\left| \int_{B_R}\left\langle \int_0^1\left( |\nabla u|^{p-2}\nabla u\right)(x+th)\,dt,\left[ |\delta_hu|^{\gamma-2}\nabla(\delta_hu)\otimes\nabla\eta \right] \cdot h  \right\rangle \,dx\right| \nonumber\\
		&\le \frac{|h|}{|h|^\gamma}\int_{B_R} \int_0^1|\nabla u(x+th)|^{p-1}\,dt\times|\delta_hu|^{\gamma-1}\left( |\nabla\eta|^2+\eta|D^2\eta|\right)\,dx\nonumber\\
		&\quad+ \frac{(\gamma-1)|h|}{|h|^\gamma}\int_{B_R} \int_0^1|\nabla u(x+th)|^{p-1}\,dt\times|\delta_hu|^{\gamma-2}\eta|\nabla(\delta_hu)||\nabla\eta|\,dx\nonumber\\
		&\le \frac{|h|}{|h|^\gamma}\|\nabla u\|^{p-1}_{L^\infty(B_R)}\int_{B_R}|\delta_hu|^{\gamma-1}\left( |\nabla\eta|^2+\eta|D^2\eta|\right)\,dx\nonumber\\
		&\quad+\frac{(\gamma-1)\|\nabla u\|^{p-1}_{L^\infty(B_R)}|h|}{|h|^\gamma}\int_{B_R}\eta|\delta_hu|^{\gamma-2}|\nabla (\delta_hu)||\nabla\eta|\,dx\nonumber\\
		&\le \|\nabla u\|^{p-1}_{L^\infty(B_R)}\int_{B_R}\frac{|\delta_hu|^{\gamma-1}}{|h|^{\gamma-1}}\left( |\nabla\eta|^2+\eta|D^2\eta|\right)\,dx\nonumber\\
		&\quad+\frac{(\gamma-1)^2}{\varepsilon}\|\nabla u\|^{2(p-1)}_{L^\infty(B_R)}\int_{B_R}\frac{|\delta_hu|^{\gamma-2}}{|h|^{\gamma-2}}|\nabla\eta|^2\,dx+\varepsilon\int_{B_R}\frac{|\delta_hu|^{\gamma-2}|\nabla(\delta_hu)|^2\eta^2}{|h|^\gamma}\,dx\nonumber\\
		&\le \varepsilon\int_{B_R}\frac{|\delta_hu|^{\gamma-2}|\nabla(\delta_hu)|^2\eta^2}{|h|^\gamma}\,dx+C(n)R^{n-2}\|\nabla u\|^{\gamma+p-2}_{L^\infty(B_R)}\nonumber\\
		&\quad+\frac{(\gamma-1)^2C(n)}{\varepsilon }R^{n-2}\|\nabla u\|^{\gamma+2p-4}_{L^\infty(B_R)}.
	\end{align}
Here $a\otimes b$ with $a,b\in\mathbb{R}^n$ denotes an $n\times n$ matrix whose $(i,j)$ entry is $a_ib_j$.

\smallskip

	\noindent\textbf{Estimate of $I_{3,2}$. }To integrate by parts without boundary terms on $B_R\times B_R$, choose $\phi_i\in C^1_0(B_R;[0,1])$ such that $\phi_i\equiv1$ in $B_{R(1-\frac{1}{i})}=:B_i$, $\phi_i\equiv0$ in $\mathbb{R}^n\backslash B_R$, and $|\nabla\phi_i|\le C(n)i/R$ for $i\ge5$. Then
	\begin{align*}
		I_{3,2}&=\frac{1}{|h|^\gamma}\int_{B_R}\int_{B_R}\frac{\left(J_q(u_h(x)-u_h(y))-J_q(u(x)-u(y))\right)}{|x-y|^{n+sq}}\\
		&\qquad\qquad\qquad\times\left[ J_\gamma(\delta_hu)(x)\eta(x)+J_\gamma(\delta_hu)(y)\eta(y)\right](\eta(x)-\eta(y))\,dxdy\\
		&=\lim\limits_{i\rightarrow+\infty}\frac{1}{|h|^\gamma}\int_{B_R}\int_{B_R}\frac{\left(J_q(u_h(x)-u_h(y))-J_q(u(x)-u(y))\right)}{|x-y|^{n+sq}}\\
		&\qquad\qquad\qquad\times\left[ J_\gamma(\delta_hu)(x)\eta(x)+J_\gamma(\delta_hu)(y)\eta(y)\right](\eta(x)-\eta(y))\phi_i(x)\phi_i(y)\,dxdy\\
		&=:\lim\limits_{i\rightarrow+\infty}(M^x_i+M^y_i),
	\end{align*}
	where $M^x_i$ and $M^y_i$ are defined by
	\begin{align*}
		M^x_i&:=\frac{1}{|h|^\gamma}\int_{B_R}\int_{B_R}\frac{\left(J_q(u_h(x)-u_h(y))-J_q(u(x)-u(y))\right)}{|x-y|^{n+sq}}\\
		&\qquad\qquad\qquad\times J_\gamma(\delta_hu)(x)\eta(x)(\eta(x)-\eta(y))\phi_i(x)\phi_i(y)\,dxdy
	\end{align*}
and
	\begin{align*}
		M^y_i&:=\frac{1}{|h|^\gamma}\int_{B_R}\int_{B_R}\frac{\left(J_q(u_h(x)-u_h(y))-J_q(u(x)-u(y))\right)}{|x-y|^{n+sq}}\\
		&\qquad\qquad\qquad\times J_\gamma(\delta_hu)(y)\eta(y)(\eta(x)-\eta(y))\phi_i(x)\phi_i(y)\,dxdy.
	\end{align*}
Hence, integration by parts and Lemma \ref{Lem5.1} give
	\begin{align*}
		M^x_i&=-\frac{1}{|h|^\gamma}\int_{B_R}\int_{B_R}\int_0^1J_q(u(x+th)-u(y+th))dt\\
		&\qquad\qquad\qquad\quad\times\nabla_{(x,y)}\left[ \frac{J_\gamma(\delta_hu)(x)\eta(x)(\eta(x)-\eta(y))\phi_i(x)\phi_i(y)}{|x-y|^{n+sq}}\right]\cdot(h,h)\,dxdy.
	\end{align*}
	Since $\nabla_{(x,y)}\frac{1}{|x-y|^{n+sq}}\cdot (h,h)=0$ and $\mathrm{supp}\, \eta\cap\mathrm{supp}\,\nabla\phi_i=\emptyset$, we deduce
	\begin{align*}
		&\quad\nabla_{(x,y)}\left[ \frac{J_\gamma(\delta_hu)(x)\eta(x)(\eta(x)-\eta(y))\phi_i(x)\phi_i(y)}{|x-y|^{n+sq}}\right]\cdot(h,h)\\
		&=\frac{(\gamma-1)|\delta_hu|^{\gamma-2}(x)\eta(x)(\eta(x)-\eta(y))\phi_i(x)\phi_i(y)\nabla(\delta_hu)(x)\cdot h}{|x-y|^{n+sq}}\\
		&\quad+\frac{J_\gamma(\delta_hu)(x)(\eta(x)-\eta(y))\phi_i(x)\phi_i(y)\nabla\eta(x)\cdot h}{|x-y|^{n+sq}}\\
		&\quad+\frac{J_\gamma(\delta_hu)(x)\eta(x)\phi_i(x)\phi_i(y)(\nabla\eta(x)-\nabla\eta(y))\cdot h}{|x-y|^{n+sq}}\\
		&\quad+\frac{J_\gamma(\delta_hu)(x)\eta(x)(\eta(x)-\eta(y))\phi_i(x)\nabla\phi_i(y)\cdot h}{|x-y|^{n+sq}}.
	\end{align*}
	Thus, after simple calculations,
	\begin{align}
		\label{5.5}
		|M^x_i|&\le\frac{C(n,\gamma,R)}{|h|^{\gamma-1}}\int_{B_R}\int_{B_R}\frac{\|\nabla u\|^{q-1}_{L^\infty(B_{2R})}}{|x-y|^{n+sq-q}}|\delta_hu|^{\gamma-2}(x)|\nabla(\delta_hu)|(x)\eta(x)\phi_i(x)\phi_i(y)\,dxdy\nonumber\\
		&\quad+\frac{C(n,R)}{|h|^{\gamma-1}}\int_{B_R}\int_{B_R}\frac{\|\nabla u\|^{q-1}_{L^\infty(B_{2R})}}{|x-y|^{n+sq-q}}|\delta_hu|^{\gamma-1}(x)\phi_i(x)\phi_i(y)\,dxdy\nonumber\\
		&\quad+\frac{C(n,R)}{|h|^{\gamma-1}}\int_{B_R}\int_{B_R}\frac{\|\nabla u\|^{q-1}_{L^\infty(B_{2R})}}{|x-y|^{n+sq-q}}|\delta_hu|^{\gamma-1}(x)\eta(x)\phi_i(x)\phi_i(y)\,dxdy\nonumber\\
		&\quad+\frac{C(n,R)}{|h|^{\gamma-1}}\int_{B_R\backslash B_i}\int_{B_{\frac{3}{4}R}}\frac{\|\nabla u\|^{q-1}_{L^\infty(B_{2R})}}{|x-y|^{n+sq-q}}|\delta_hu|^{\gamma-1}i\,dxdy\nonumber\\
		&=:M_{i,1}^x+M_{i,2}^x+M_{i,3}^x+M_{i,4}^x.
	\end{align}
For $M^x_{i,1}$, by Young's inequality and noting $n+sq-q<n$, we have
	\begin{align}
		\label{5.6}
		M^x_{i,1}&\le \frac{C(n,s,q,\gamma,R)}{|h|^{\gamma-1}}\|\nabla u\|^{q-1}_{L^\infty(B_{2R})}\int_{B_R}|\delta_hu|^{\gamma-2}|\nabla(\delta_hu)|\eta\,dx\nonumber\\
		&\le \frac{\varepsilon}{2}\int_{B_R}\frac{|\delta_hu|^{\gamma-2}|\nabla(\delta_hu)|^2\eta^2}{|h|^\gamma}\,dx+\frac{C^2(n,s,q,\gamma,R)\|\nabla u\|^{2(q-1)}_{L^\infty(B_{2R})}}{\varepsilon}\int_{B_R}\frac{|\delta_hu|^{\gamma-2}}{|h|^{\gamma-2}}\,dx\nonumber\\
		&\le \frac{\varepsilon}{2}\int_{B_R}\frac{|\delta_hu|^{\gamma-2}|\nabla(\delta_hu)|^2\eta^2}{|h|^\gamma}\,dx+\frac{C^2(n,s,q,\gamma,R)\|\nabla u\|^{\gamma+2q-4}_{L^\infty(B_{2R})}}{\varepsilon}.
	\end{align}
For $M^x_{i,2}$ and $M^x_{i,3}$, the Lipschitz continuity of $u$ implies
	\begin{align}
		\label{5.7}
		M^x_{i,2}+M^x_{i,3}\le C(n,s,q,R)\|\nabla u\|^{\gamma+q-2}_{L^\infty(B_{2R})}.
	\end{align}
For $M^x_{i,4}$, note that $x\in B_{\frac{3}{4}R}$ and $y\in B_R\backslash B_i$ imply $|x-y|\ge|y|-|x|\ge \frac{4}{5}R-\frac{3}{4}R=  \frac{R}{20}$. Thus,
	\begin{align}
		\label{5.8}
		M^x_{i,4}&\le C\|\nabla u\|^{\gamma+q-2}_{L^\infty(B_{2R})}\int_{B_R\backslash B_i}\int_{B_{\frac{3}{4}R}}i\,dxdy\nonumber\\
		&\le C\|\nabla u\|^{\gamma+q-2}_{L^\infty(B_{2R})}\frac{1-\left( 1-\frac{1}{i}\right)^n }{\frac{1}{i}}\le C\|\nabla u\|^{\gamma+q-2}_{L^\infty(B_{2R})},
	\end{align}
		where we used $\lim\limits_{i\rightarrow+\infty}\frac{1-\left( 1-\frac{1}{i}\right)^n }{\frac{1}{i}}=n$, and $C>0$ depends on $n,s,q,R$.
	
Substituting \eqref{5.6}--\eqref{5.8} into \eqref{5.5} arrives at
	\begin{align*}
		|M^x_i|&\le \frac{\varepsilon}{2} \int_{B_R}\frac{|\delta_hu|^{\gamma-2}|\nabla(\delta_hu)|^2\eta^2}{|h|^\gamma}\,dx+\frac{C}{\varepsilon}\|\nabla u\|^{\gamma+2q-4}_{L^\infty(B_{2R})}+C\|\nabla u\|^{\gamma+q-2}_{L^\infty(B_{2R})},
	\end{align*}
	where $C>0$ depends on $n,s,q,\gamma,R$. Similarly, by Fubini's theorem, we also have
	\begin{align*}
		|M^y_i|&\le \frac{\varepsilon}{2} \int_{B_R}\frac{|\delta_hu|^{\gamma-2}|\nabla(\delta_hu)|^2\eta^2}{|h|^\gamma}\,dx+\frac{C}{\varepsilon}\|\nabla u\|^{\gamma+2q-4}_{L^\infty(B_{2R})}+C\|\nabla u\|^{\gamma+q-2}_{L^\infty(B_{2R})}.
	\end{align*}
Therefore,
	\begin{align}
		\label{6.1}
		|I_{3,2}|&\le \varepsilon \int_{B_R}\frac{|\delta_hu|^{\gamma-2}|\nabla(\delta_hu)|^2\eta^2}{|h|^\gamma}\,dx+\frac{C}{\varepsilon}\|\nabla u\|^{\gamma+2q-4}_{L^\infty(B_{2R})}+C\|\nabla u\|^{\gamma+q-2}_{L^\infty(B_{2R})}.
	\end{align}
Here the constant $C>0$ depends on $n,s,q,\gamma,R$.

\smallskip

	\noindent\textbf{Estimate of $I_{4}$. }To avoid boundary terms in the integration by parts, choose $\phi^i\in C^1_0(\mathbb{R}^n;[0,1])$, $i\ge1$, such that $\phi^i\equiv1$ in $B_{iR}\backslash B_{R(1+1/i)}$, $\phi^i\equiv0$ in $(\mathbb{R}^n\backslash B_{2iR})\cup B_R$, $|\nabla\phi^i|\le C(n)i/R$ in $B_{R(1+1/i)}\backslash B_R$, and $|\nabla\phi^i|\le C(n)/(iR)$ in $B_{2iR}\backslash B_{iR}$. Then
	\begin{align*}
		I_4&=\frac{1}{|h|^\gamma}\int_{\mathbb{R}^n\backslash B_R}\int_{B_{\frac{3}{4}R}}\frac{\left(J_q(u_h(x)-u_h(y))-J_q(u(x)-u(y))\right)\eta^2(x)J_\gamma(\delta_hu)(x)}{|x-y|^{n+sq}}\,dxdy\\
		&=\lim\limits_{i\rightarrow+\infty}\int_{\mathbb{R}^n\backslash B_R}\int_{B_{\frac{3}{4}R}}\frac{\left(J_q(u_h(x)-u_h(y))-J_q(u(x)-u(y))\right)\eta^2(x)J_\gamma(\delta_hu)(x)\phi^i(y)}{|h|^\gamma|x-y|^{n+sq}}\,dxdy\\
		&=:\lim\limits_{i\rightarrow+\infty}T_i.
	\end{align*}
	By integration by parts, we have
	\begin{align}
		\label{5.90}
		T_i&=-\int_{\mathbb{R}^n\backslash B_R}\int_{B_{\frac{3}{4}R}}\int_0^1J_q(u(x+th)-u(y+th))\,dt\nonumber\\
		&\qquad\qquad\qquad\quad\times\nabla_{(x,y)}\left[ \frac{\eta^2(x)J_\gamma(\delta_hu)(x)\phi^i(y)}{|x-y|^{n+sq}}\right]\cdot(h,h)\,dxdy,
	\end{align}
		where, for $x\neq y$, we have
	\begin{align*}
		&\quad\nabla_{(x,y)}\left[ \frac{\eta^2(x)J_\gamma(\delta_hu)(x)\phi^i(y)}{|x-y|^{n+sq}}\right]\cdot(h,h)\\
		&=\frac{2\eta(x)J_\gamma(\delta_hu)(x)\phi^i(y)\nabla\eta(x)\cdot h}{|x-y|^{n+sq}}+\frac{(\gamma-1)\eta^2(x)|\delta_hu|^{\gamma-2}(x)\phi^i(y)\nabla(\delta_hu)(x)\cdot h}{|x-y|^{n+sq}}\\
		&\quad+\frac{\eta^2(x)J_\gamma(\delta_hu)(x)\nabla\phi^i(y)\cdot h}{|x-y|^{n+sq}}.
	\end{align*}
	Substituting this expression into \eqref{5.90} and noting that $|x+th|\le \left( \frac{3}{4}+\frac{1}{16}\right)R $ and $|y+th|\ge \left( 1-\frac{1}{16}\right)R $, we obtain
	\begin{align}
		\label{5.9}
		|T_i|&\le \frac{2}{|h|^{\gamma-1}}\int_{\mathbb{R}^n\backslash B_{\frac{15}{16}R}}\int_{B_{\frac{13}{16}R}}\frac{|J_q(u(x)-u(y))|}{|x-y|^{n+sq}}\eta(x)|\nabla\eta(x)||\delta_hu|^{\gamma-1}(x)\,dxdy\nonumber\\
		&\quad+ \frac{\gamma-1}{|h|^{\gamma-1}}\int_{{\mathbb{R}^n\backslash B_{\frac{15}{16}R}}}\int_{B_{\frac{13}{16}R}}\frac{|J_q(u(x)-u(y))|}{|x-y|^{n+sq}}\eta^2(x)|\delta_hu|^{\gamma-2}(x)|\nabla(\delta_hu)|(x)\,dxdy\nonumber\\
		&\quad+ \frac{1}{|h|^{\gamma-1}}\int_{{\mathbb{R}^n\backslash B_{\frac{15}{16}R}}}\int_{B_{\frac{13}{16}R}}\frac{|J_q(u(x)-u(y))|}{|x-y|^{n+sq}}\eta^2(x)|\delta_hu|^{\gamma-1}(x)|\nabla \phi^i(y)|\,dxdy\nonumber\\
		&=:T_{i,1}+T_{i,2}+T_{i,3}.
	\end{align}
	
	For $T_{i,1}$, observe that $|x-y|\ge |y|-\frac{13}{16}R\ge |y|-\frac{13}{15}|y| =\frac{2}{15}|y|$ for $x\in B_{\frac{13}{16}R}$. Hence,
	\begin{align}
		\label{T}
		\int_{\mathbb{R}^n\backslash B_{\frac{15}{16}R}}\frac{|u(x)-u(y)|^{q-1}}{|x-y|^{n+sq}}\,dy&
		\le \int_{\mathbb{R}^n\backslash B_{\frac{15}{16}R}}\left( \frac{\|u\|^{q-1}_{L^\infty(B_R)}}{|x-y|^{n+sq}}+\frac{|u(y)|^{q-1}}{|x-y|^{n+sq}}\right)\,dy\nonumber\\
		&\le C\|u\|^{q-1}_{L^\infty(B_R)}+C\int_{\mathbb{R}^n\backslash B_{\frac{15}{16}R}} \frac{|u(y)|^{q-1}}{|y|^{n+sq}}\,dy\nonumber\\
		&\le C\|u\|^{q-1}_{L^\infty(B_R)}+C\mathrm{Tail}\left( u;x_0,\frac{15}{16}R\right) ^{q-1}
		\le C\mathfrak{M}_{2R}^{q-1},
	\end{align}
	 where $C>0$ depends on $n,q,s,R$ and in the last line we used the following inequality
	\begin{align*}
		\mathrm{Tail}\left( u;x_0,\frac{15}{16}R\right) ^{q-1}\le C(n)\left(\|u\|_{L^\infty(B_{2R})}+\mathrm{Tail}(u;x_0,2R)\right)^{q-1}.
	\end{align*}
	As a result,
	\begin{align}
		\label{5.10}
		T_{i,1}&\le C(n,s,q,R)\mathfrak{M}_{2R}^{q-1}\int_{B_{\frac{13}{16}R}}\frac{|\delta_hu|^{\gamma-1}}{|h|^{\gamma-1}}\,dx\le C(n,s,q,R)\mathfrak{M}_{2R}^{q-1}\|\nabla u\|^{\gamma-1}_{L^\infty(B_R)}.
	\end{align}
	
	Similar to \eqref{T}, by Young's inequality and the local Lipschitz continuity of $u$, for $T_{i,2}$ we infer that
	\begin{align}
		\label{5.11}
		T_{i,2}&\le C(n,s,q,\gamma,R)\mathfrak{M}_{2R}^{q-1}\int_{B_R}\frac{|\delta_hu|^{\gamma-2}|\nabla(\delta_hu)|\eta^2}{|h|^{\gamma-1}}\,dx \nonumber\\
		&\le \varepsilon \int_{B_R}\frac{|\delta_hu|^{\gamma-2}|\nabla(\delta_hu)|^2\eta^2}{|h|^\gamma}\,dx+\frac{C^2(n,s,q,\gamma,R)}{\varepsilon}\mathfrak{M}_{2R}^{2q-2}\|\nabla u\|^{\gamma-2}_{L^\infty(B_R)}.
	\end{align}
	
	Observe that $\lim\limits_{i\rightarrow+\infty}i|B_{R\left( 1+\frac{1}{i}\right) }\backslash B_R|=nR^n$. For $T_{i,3}$, we have
	\begin{align}
		\label{5.12}
		T_{i,3}&\le \frac{C(n)}{|h|^{\gamma-1}}\int_{B_{2iR}\backslash B_{iR}}\int_{B_{\frac{13}{16}R}}\frac{1}{iR}\frac{|u(x)-u(y)|^{q-1}}{|x-y|^{n+sq}}\eta^2(x)|\delta_hu|^{\gamma-1}(x)\,dxdy\nonumber\\
		&\quad+\frac{C(n)}{|h|^{\gamma-1}}\int_{B_{R\left( 1+\frac{1}{i}\right) }\backslash B_{R}}\int_{B_{\frac{13}{16}R}}\frac{i}{R}\frac{|u(x)-u(y)|^{q-1}}{|x-y|^{n+sq}}\eta^2(x)|\delta_hu|^{\gamma-1}(x)\,dxdy\nonumber\\
		&\le \frac{C(n,s,q,R)}{|h|^{\gamma-1}}\int_{\mathbb{R}^n\backslash B_R}\int_{B_{\frac{13}{16}R}}\frac{1}{i}\frac{|u(x)-u(y)|^{q-1}}{|x-y|^{n+sq}}\eta^2(x)|\delta_hu|^{\gamma-1}(x)\,dxdy\nonumber\\
		&\quad+C(n,s,q,R)\int_{B_{R\left( 1+\frac{1}{i}\right) }\backslash B_{R}}i\|u\|^{q-1}_{L^\infty(B_{2R})}\|\nabla u\|^{\gamma-1}_{L^\infty(B_{R})}\,dy\nonumber\\
		&\le C(n,s,q,R)\mathfrak{M}_{2R}^{q-1}\|\nabla u\|^{\gamma-1}_{L^\infty(B_{R})},
	\end{align}
by letting $i\rightarrow\infty$. Here we need to notice the integral in the third line will tend to 0 as $i\rightarrow\infty$. Consequently, we substitute \eqref{5.10}--\eqref{5.12} into \eqref{5.9} to get
	\begin{align*}
		|T_i|&\le \varepsilon \int_{B_R}\frac{|\delta_hu|^{\gamma-2}|\nabla(\delta_hu)|^2\eta^2}{|h|^\gamma}\,dx+\frac{C}{\varepsilon}\|\nabla u\|^{\gamma-2}_{L^\infty(B_{R})}\mathfrak{M}_{2R}^{2q-2}
		+ C\|\nabla u\|^{\gamma-1}_{L^\infty(B_{R})}\mathfrak{M}_{2R}^{q-1}
	\end{align*}
	with $C>0$ depending on $n,s,q,\gamma,R$. Therefore, one can estimate $I_4$ as
	\begin{align}
		\label{6.2}
		|I_4|&= \lim\limits_{i\rightarrow+\infty}|T_i|\nonumber\\
		&\le \varepsilon \int_{B_R}\frac{|\delta_hu|^{\gamma-2}|\nabla(\delta_hu)|^2\eta^2}{|h|^\gamma}\,dx+\frac{C}{\varepsilon}\|\nabla u\|^{\gamma-2}_{L^\infty(B_{R})}\mathfrak{M}_{2R}^{2q-2}+ C\|\nabla u\|^{\gamma-1}_{L^\infty(B_{R})}\mathfrak{M}_{2R}^{q-1}
	\end{align}
	with $C>0$ depending on $n,s,q,\gamma,R$.

	Combining \eqref{5.3}, \eqref{5.4}, \eqref{6.1}, and \eqref{6.2} gives
	\begin{align*}
		I_1&=\frac{1}{|h|^\gamma}\int_{B_R}\eta^2|\delta_hu|^{\gamma-2}\left(  |\nabla u_h|^{p-2}\nabla u_h-|\nabla u|^{p-2}\nabla u\right)\cdot(\nabla u_h-\nabla u)\,dx\\
		&\le 3\varepsilon \int_{B_R}\frac{|\delta_hu|^{\gamma-2}|\nabla(\delta_hu)|^2\eta^2}{|h|^\gamma}\,dx+\frac{C}{\varepsilon}\|\nabla u\|^{\gamma-2}_{L^\infty(B_{2R})}\left( \mathfrak{M}_{2R}^{2q-2}+\|\nabla u\|^{2q-2}_{L^\infty(B_{2R})}\right) \\
		&\quad+C\|\nabla u\|^{\gamma-1}_{L^\infty(B_{2R})}\left( \mathfrak{M}_{2R}^{q-1}+\|\nabla u\|^{q-1}_{L^\infty(B_{2R})}\right).
	\end{align*}
	 Here $C>0$ depends on $n,s,q,\gamma,R$. Next, we estimate $I_1$ from below. We utilize the following inequality
	  \begin{align*}
	  	C(p)|\nabla u_h-\nabla u|^2
	  	\le \left( 1+|\nabla u_h|^2+|\nabla u|^2\right)^{1-\frac{p}{2}} \left(  |\nabla u_h|^{p-2}\nabla u_h-|\nabla u|^{p-2}\nabla u\right)\cdot(\nabla u_h-\nabla u)
	  \end{align*}
	  to obtain
	  \begin{align}
	  	\label{3.07}
	  	&\quad\frac{(p-1)}{|h|^\gamma}\int_{B_R}\eta^2|\nabla u_h-\nabla u|^2|u_h-u|^{\gamma-2}\,dx\nonumber\\
	  	&\le \frac{1+\|\nabla u\|^{2-p}_{L^\infty(B_{2R})}}{|h|^\gamma} \int_{B_R}\eta^2|u_h-u|^{\gamma-2}\left(  |\nabla u_h|^{p-2}\nabla u_h-|\nabla u|^{p-2}\nabla u\right)\cdot(\nabla u_h-\nabla u)\,dx\nonumber\\
	  	&=\left( 1+\|\nabla u\|^{2-p}_{L^\infty(B_{2R})}\right)I_1.
	  \end{align}
	 Thus, we find
	  \begin{align*}
	  	&\quad\frac{(p-1)}{|h|^\gamma}\int_{B_R}\eta^2|\nabla u_h-\nabla u|^2|u_h-u|^{\gamma-2}\,dx\\
	  	&\le 3\varepsilon \left( 1+\|\nabla u\|^{2-p}_{L^\infty(B_R)}\right)\int_{B_R}\frac{|\delta_hu|^{\gamma-2}|\nabla(\delta_hu)|^2\eta^2}{|h|^\gamma}\,dx\\
	  	&\quad+\frac{C(n,s,q,\gamma,R)}{\varepsilon}\|\nabla u\|^{\gamma-2}_{L^\infty(B_{2R})}\left( 1+\|\nabla u\|^{2-p}_{L^\infty(B_{2R})}\right)\left( \mathfrak{M}_{2R}^{2q-2}+\|\nabla u\|^{2q-2}_{L^\infty(B_{2R})}\right) \\
	  	&\quad+C(n,s,q,R)\|\nabla u\|^{\gamma-1}_{L^\infty(B_{2R})}\left( 1+\|\nabla u\|^{2-p}_{L^\infty(B_{2R})}\right)\left( \mathfrak{M}_{2R}^{q-1}+\|\nabla u\|^{q-1}_{L^\infty(B_{2R})}\right).
	  \end{align*}
	  By setting
	  \begin{align*}
	  	\varepsilon=\frac{p-1}{6\left( 1+\|\nabla u\|^{2-p}_{L^\infty(B_{2R})}\right)},
	  \end{align*}
	  one has
	  \begin{align*}
	  	&\quad\frac{(p-1)}{2|h|^\gamma}\int_{B_R}\eta^2|\nabla u_h-\nabla u|^2|u_h-u|^{\gamma-2}\,dx\\
	  	&\le C(n,s,q,\gamma,R)\|\nabla u\|^{\gamma-2}_{L^\infty(B_{2R})}\left( 1+\|\nabla u\|^{2-p}_{L^\infty(B_{2R})}\right)^2\left( \mathfrak{M}_{2R}^{2q-2}+\|\nabla u\|^{2q-2}_{L^\infty(B_{2R})}\right) \\
	  	&\quad+C(n,s,q,R)\|\nabla u\|^{\gamma-1}_{L^\infty(B_{2R})}\left( 1+\|\nabla u\|^{2-p}_{L^\infty(B_{2R})}\right)\left( \mathfrak{M}_{2R}^{q-1}+\|\nabla u\|^{q-1}_{L^\infty(B_{2R})}\right).
	  \end{align*}
	Note that $\nabla J_{\frac{\gamma}{2}+1}(\delta_h u)=\frac{\gamma}{2}|\delta_hu|^\frac{\gamma-2}{2}\nabla \delta_hu$. The preceding display becomes
	  \begin{align}
	  	\label{eq:subquadratic-gradient-energy}
	  	\int_{B_R}\eta^2\left| \nabla J_{\frac{\gamma}{2}+1}(u_h-u)\right|^2\,dx\le C(n,s,p,q,\gamma,R)\left( 1+\mathfrak{M}_{2R}+\|\nabla u\|_{L^\infty(B_{2R})}\right)^{\gamma-2p+2q} |h|^{\gamma}.
	  \end{align}
	
	  Applying Lemma \ref{Lem3.7}, and observing that $J_{\frac{\gamma}{2}+1}(\delta_hu)\in W^{1,2}(B_R)$ and $|\delta_hJ_{\frac{\gamma}{2}+1}(\delta_hu)|\ge |\delta_h(\delta_hu)|^{\frac{\gamma}{2}}$, we obtain
	  \begin{align*}
	  	\int_{B_\frac{R}{4}}|\delta_h(\delta_hu)|^\gamma\,dx&\le |h|^2\int_{B_\frac{R}{2}}\left| \nabla J_{\frac{\gamma}{2}+1}(u_h-u)\right|^2\,dx\\
	  	&\le |h|^2\int_{B_R}\eta^2\left| \nabla J_{\frac{\gamma}{2}+1}(u_h-u)\right|^2\,dx\\
	  	&\le C(n,s,p,q,\gamma,R)\left( 1+\mathfrak{M}_{2R}+\|\nabla u\|_{L^\infty(B_{2R})}\right)^{\gamma-2p+2q} |h|^{\gamma+2}.
	  \end{align*}
	  For any $B_{2R}(x_0)\subset\subset\Omega$, Lemma \ref{Lem2.18} yields $\nabla u\in W^{\beta,\gamma}(B_{\frac{R}{8}})$
	  and furthermore,
	  \begin{align*}
	  	[\nabla u]^\gamma_{W^{\beta,\gamma}(B_{\frac{R}{8}})}\le C(n,s,p,q,\gamma,\beta,R)\left( 1+\mathfrak{M}_{2R}+\|\nabla u\|_{L^\infty(B_{2R})}\right)^{\gamma-2p+2q}
	  \end{align*}
	  for any $\beta\in\left( 0,\frac{2}{\gamma}\right) $. Thus, we conclude $\nabla u\in W^{\beta,\gamma}_\mathrm{loc}(\Omega)$
	  for any $\gamma\ge2$ and $\beta\in \left( 0,\frac{2}{\gamma}\right)$.
\end{proof}

Taking $\gamma=2$ in the preceding argument yields the $W^{2,2}_{\rm loc}$ regularity; in this case, second-order differences are unnecessary.

\smallskip

\begin{proof}[Proof of Theorem \ref{th1} \rm (ii)]
	In the case $\gamma=2$, there is no need to consider the relationship between $\nabla J_{\frac{\gamma}{2}+1}(\delta_hu)$ and $\delta_hJ_{\frac{\gamma}{2}+1}(\delta_hu)$. Therefore, one can readily infer from \eqref{eq:subquadratic-gradient-energy} that
	\begin{align*}
		\int_{B_{\frac{R}{2}}(x_0)}\left| \nabla u_h-\nabla u\right|^2\,dx&\le\int_{B_R(x_0)}\eta^2\left| \nabla u_h-\nabla u\right|^2\,dx\\
		&\le C(n,s,p,q,\gamma,R)\left( 1+\mathfrak{M}_{2R}+\|\nabla u\|_{L^\infty(B_{2R})}\right)^{2-2p+2q} |h|^2.
	\end{align*}
	The preceding inequality and the relation between differences and derivatives imply $\nabla u\in W^{1,2}\left(B_\frac{R}{8}(x_0);\mathbb{R}^n\right)$, and
	\begin{align*}
		[\nabla u]^2_{W^{1,2}(B_{\frac{R}{8}}(x_0))}\le C(n,s,p,q,\gamma,R)\left( 1+\mathfrak{M}_{2R}+\|\nabla u\|_{L^\infty(B_{2R})}\right)^{2-2p+2q}.
	\end{align*}
Hence $u\in W^{2,2}_\mathrm{loc}(\Omega)$, which completes the proof.
\end{proof}

\smallskip

\subsection*{Acknowledgments}
This work was supported by the National Natural Science Foundation of China (Nos. 12471128 and 12301245) and the Natural Science Foundation of Heilongjiang Province (No. YQ2025A002).

\subsection*{Conflict of interest}  The authors declare no conflict of interest.
	
\subsection*{Data availability}
No data were generated or analyzed in this study.

\end{document}